\documentclass{amsart}

\usepackage{ifthen}
\usepackage{verbatim}

\usepackage[curve]{xypic}
  \ifthenelse{\isundefined{\xymatrix}}%
     {\newif\ifnoxypic\noxypictrue}%
     {\newif\ifnoxypic\noxypicfalse}
\def\xypicmessage{\centerline{\fbox{\color{red}{\LARGE XY-Pic
package needs to by loaded.}}}}

\usepackage{amssymb}
\usepackage{amsmath}
\usepackage{amsthm}
\usepackage{amsfonts}
\usepackage{amscd}
\usepackage{mathtools}
\usepackage{boxedminipage}
\usepackage[usenames]{color}
\usepackage{graphicx} 
\usepackage{colordvi}
\newif\ifhrule\hrulefalse
\newtheorem{Lemma}{Lemma}[section]
\newtheorem{Theorem}[Lemma]{Theorem}
\newtheorem{Corollary}[Lemma]{Corollary}
\newtheorem{Proposition}[Lemma]{Proposition}
\theoremstyle{definition}
\newtheorem{Definition}[Lemma]{Definition}
\newtheorem{Remark}[Lemma]{Remark}

\numberwithin{equation}{section}

\newcommand{\mb}{\mathbb}
\newcommand{\mc}{\mathcal}

\newcommand{\bs}{\boldsymbol}
\newcommand{\CoxeterGroup}[1]{\mc #1}
\newcommand{\Sn}[1][n]{\CoxeterGroup{S_{#1}}}
\newcommand{\Bn}[1][n]{\CoxeterGroup{B_{#1}}}
\newcommand{\Dn}[1][n]{\CoxeterGroup{D_{#1}}}
\newcommand{\Arrangement}[1][A]{\mc{#1}}

\newcommand{\Faces}{\mc F}

\newcommand{\DescentAlgebra}[1][W]{\ensuremath{\Sigma_k(#1)}}

\newcommand{\InvariantSubalgebra}[1][W]{(k\Faces)^{#1}}

\newcommand{\InvariantSubalgebraB}{\InvariantSubalgebra[\Bn]}

\newcommand{\InvariantQuiver}{\Gamma}
\newcommand{\InvariantQuiverW}[1][W]{\InvariantQuiver_{#1}}

\newcommand{\InvariantQuiverD}[1][n]{\InvariantQuiverW[\Dn]}
\newcommand{\IntersectionLattice}{\mc L}
\newcommand{\Orbit}{\mc O}
\newcommand{\Quiver}{\mc Q}

\newcommand{\SetPartitionB}{\ensuremath{\{B_1,\ldots,B_r;Z\}}}

\newcommand{\Partition}[1]{\pi(X_{#1})}
\newcommand{\PartitionLattice}[1][A]{\ensuremath{\Pi^{#1}_n}}

\newcommand{\PartitionLatticeB}{\PartitionLattice[B]}

\newcommand{\supp}{\operatorname{supp}}

\newcommand{\rank}{\operatorname{rank}}

\newcommand{\End}{\operatorname{End}}

\newcommand{\rad}{\operatorname{rad}}

\newcommand{\img}{\operatorname{im}}

\newcommand{\inv}{^{-1}}
\newcommand{\opp}{^\ast}
\newcommand{\scissors}{\ensuremath{^\wedge}}

\newcommand{\CompleteSystemPunctuated}[1][c]{{#1}omplete system of primitive orthogonal idempotents}
\newcommand{\CompleteSystem}[1][c]{{\CompleteSystemPunctuated[#1] }}

\newcommand{\RAD}[1][]{\ensuremath{\rad^{#1}\InvariantSubalgebra[\Sn]}}
\newcommand{\Norm}{{\mc{N}}}
\newcommand{\iNorm}[1]{\varphi(\Norm(#1))}
\renewcommand{\eqref}[1]{(\ref{#1})}
\newcommand{\InvariantQuiverAlgebra}[1][W]{\ensuremath{(k\Quiver)^{#1}}}

\newenvironment{note}[1][Note]
 {\begin{center}\begin{boxedminipage}{4.5in}\setlength{\parindent}{1em}\noindent\textbf{#1. }}
 {\end{boxedminipage}\end{center}}

\newcommand{\cldots}{\mathord{.}\mathord{.}\mathord{.}}
\newcommand{\ccdots}{\mathord{\cdot}\mathord{\cdot}\mathord{\cdot}}

\ifthenelse{\isundefined{\arrow}}%
 {\newcommand{\arrow}{\mathord{\rightarrow}}}%
 {\renewcommand{\arrow}{\mathord{\rightarrow}}}%
\newcommand{\labelledarrow}[1]{\mathord{\buildrel{#1}\over{\rightarrow}}}
\newcommand{\arrows}{\arrow\ccdots\arrow}
\newcommand{\Path}[3][1]{
\ifthenelse{\equal{#1}{2}}%
 {\ensuremath{\left({#2}_0\arrow{#2}_1\arrows{#2}_{#3}\right)}}%
 {\ensuremath{\left({#2}_0\arrows{#2}_{#3}\right)}}%
}
\newcommand{\adots}{\mathord{\cdot}\mathord{\cdot}\mathord{\cdot}}

\newcommand{\defn}[1]{\textit{\textbf{#1}}}

\allowdisplaybreaks
\begin{document}

%%%%%
% Some notes and warnings.

% Warns that XY-Pic is not loaded.
\ifnoxypic\xypicmessage\vspace{0.5in}\fi

%%%%%
% Title

\title{On the quiver of the descent algebra}
\author{Franco V. Saliola}
\address{%
	Laboratoire de Combinatoire et d'Informatique Math\'ematique\\
	Universit\'e du Qu\'ebec \`a Montr\'eal\\
	Case Postale 8888, succursale Centre-ville\\
	Montr\'eal, Qu\'ebec H3C 3P8\\
	Canada}
\thanks{This research was supported, in part, by an NSERC PGS B grant.}
\email{saliola@gmail.com}
\dedicatory{Dedicated to the memory of Manfred Schocker (1970-2006).}

\begin{abstract}
We study the quiver of the descent algebra of a finite Coxeter group
$W$. The results include a derivation of the quiver of the descent
algebra of types $A$ and $B$. Our approach is to study the descent
algebra as an algebra constructed from the reflection arrangement
associated to $W$. 
\end{abstract}

\maketitle

%%%%%
% Header information.
\markboth{\textsc{Franco V. Saliola}}
         {\textsc{On the quiver of the descent algebra}}

\tableofcontents

\section{Introduction}

The descent algebra $\DescentAlgebra$ of a finite Coxeter group $W$ is
a highly exceptional subalgebra of the group algebra of $W$. First
introduced by Louis Solomon in 1976 \cite{Solomon1976}, it has enjoyed
much attention because of several connections with various areas of
mathematics, including the representation theory of Coxeter groups,
free Lie algebras and higher Lie modules, Hochschild homology, and
probability. These connections are described in a survey article by
Manfred Schocker \cite{Schocker2004}.

In this article we study the quiver of the descent algebra. Our
approach is to use a result of T.\ P.\ Bidigare that identities the
descent algebra with the $W$-invariant subalgebra
$\InvariantSubalgebra$ of a semigroup algebra $k\Faces$ associated to
the reflection arrangement of $W$ \cite{Bidigare1997}. Then using
results about $k\Faces$ we deduce some general properties about the
quiver of the descent algebra and determine the quiver of the descent
algebras of type $A$ and $B$. The quiver of the descent algebra of
type $A$ has already been computed \cite{Schocker2004}, but the quiver
of the descent algebra of type $B$ was not previously known.

We briefly outline the contents and structure of the article. Section
\ref{s:GeometricApproach} defines finite Coxeter groups and reflection
arrangements, and explains the connection between descent algebra
$\DescentAlgebra$ and the $W$-invariant subalgebra
$\InvariantSubalgebra$. Section \ref{s:QuiverOfSplitBasicAlgebra}
recalls definitions and results about quivers of split basic algebras.
In Section \ref{s:MyAlgebrasAreBasic} we provide a proof that
$k\Faces$ and $\InvariantSubalgebra$ are split basic algebras, so
there are (canonical) quivers associated to each. Section
\ref{s:CompleteSystems} constructs a \CompleteSystem for $k\Faces$
that leads to a \CompleteSystem for $\InvariantSubalgebra$. This
allows us to define, in Section \ref{s:WEquivariantSurjection}, a
$W$-equivariant surjection $\varphi:k\Quiver\twoheadrightarrow
k\Faces$, where $\Quiver$ is the quiver of $k\Faces$. We use this
surjection in Section \ref{s:OnTheQuiverOfTheInvarviantSubalgebra} to
prove some general properties of the quiver of the descent algebra and
in Sections \ref{s:QuiverDesAlgA} and \ref{s:QuiverDesAlgB}   to
determine the quiver of the descent algebras of types $A$ and $B$,
respectively. Finally, Section \ref{s:FutureDirections} discusses some
future directions for this project.

The interested reader may also decide to consult recent work of G\"otz
Pfeiffer who is taking a different approach to the problem of
determining the quiver of the descent algebras \cite{Pfeiffer2007}. 

\section{The Geometric Approach to the Descent Algebra}
\label{s:GeometricApproach}

For an introduction to the theory of Coxeter groups, see the books
\cite{Brown1989,Humphreys1990,Kane2001,BjornerBrenti2005}. The reader
may wish to read \S\ref{ss:CoxSys+RefArrs} and
\S\ref{ss:geometricApproach} alongside
\S\ref{ss:SymmetricGroupExample} since the latter presents these ideas
for the symmetric group $\Sn$. Also see \S\ref{ss:TypeBCoxGrp}, which
describes some of these ideas for the hyperoctahedral group $\Bn$.

\subsection{Coxeter systems and reflection arrangements}
\label{ss:CoxSys+RefArrs}

Let $V$ be a finite dimensional real vector space. A \defn{finite
Coxeter group} $W$ is a finite group generated by a set of reflections
of $V$. The \defn{reflection arrangement} of $W$ is the hyperplane
arrangement $\Arrangement$ consisting of the hyperplanes of $V$ fixed
by some reflection in $W$. The connected components of the complement
of $\bigcup_{H\in\Arrangement} H$ in $V$ are called \defn{chambers}. A
\defn{wall} of a chamber $c$ is a hyperplane $H\in\Arrangement$ such
that $H \cap \overline c$ spans $H$, where $\overline c$ denotes the
closure of the set $c$. 

Fix a chamber $c$ and let $S \subseteq W$
denote the set of reflections in the walls of $c$. Then $S$ is a
generating set of $W$ \cite[{\S}I.5A]{Brown1989}. The pair $(W,S)$ is
called a \defn{Coxeter system}, and $c$ is the \defn{fundamental
chamber} of $(W,S)$.

\subsection{The descent algebra}
\label{ss:DescentAlgebra}

Fix a Coxeter system $(W,S)$. For $J \subseteq S$, let $W_J = \langle
J \rangle$ denote the subgroup of $W$ generated by the elements in
$J$. Each coset of $W_J$ in $W$ contains a unique element of minimal
length, where the \defn{length} $\ell(w)$ of an element $w$ of $W$ is
the smallest number of generators $s_1, \ldots, s_i \in S$ such that
$w = s_1 \cdots s_i$ \cite[Proposition 1.10(c)]{Humphreys1990}. 

For $J\subseteq S$, let $X_J$ denote the set of \defn{minimal length
coset representatives} of $W_J$ and let $x_J = \sum_{w \in X_J} w$
denote the sum of the elements of $X_J$ in the group algebra $kW$ of
$W$ with coefficients in a field $k$. Louis Solomon proved that the
elements $x_J$ form a $k$-vector space basis of a subalgebra of $kW$
\cite[Theorem 1]{Solomon1976}. This subalgebra is denoted by
$\DescentAlgebra$ and is called the \defn{descent algebra} of $W$.
Throughout $k$ will be a field of characteristic that does not divide
the order of $W$.

\subsection{The geometric approach to the descent algebra}
\label{ss:geometricApproach}

Let $(W,S)$ be a finite Coxeter system with fundamental chamber $c$,
and let $\Arrangement$ be the reflection arrangement of $W$. The
reader may want to read this section alongside Section
\ref{ss:SymmetricGroupExample}.

\subsubsection{Face semigroup algebra}
\label{sss:FaceSemigroupAlgebra}
For each hyperplane $H \in \Arrangement$, let $H^+$ and $H^-$ denote
the two open half spaces of $V$ determined by $H$. The choice of
labels $H^+$ and $H^-$ is arbitrary, but fixed throughout. For
convenience, let $H^0=H$. A \defn{face} of $\Arrangement$ is a
non-empty intersection of the form $\bigcap_{H \in \Arrangement}
H^{\sigma_H}$, where $\sigma_H \in \{+,0,-\}$ for each hyperplane $H
\in \Arrangement$. The sequence $(\sigma_H)_{H \in \Arrangement}$ is
called the \defn{sign sequence} of the face. We denote the sign
sequence of a face $x$ by $\sigma(x)=(\sigma_H(x))_{H \in \Arrangement}$. 

Let $\Faces$ denote the set of all faces of $\Arrangement$. 
Define the product of two faces $x, y \in \Faces$ to be the face $xy$
with sign sequence 
\begin{gather*}
\sigma_H(xy) = 
\begin{cases}
\sigma_H(x), & \text{if } \sigma_H(x) \neq 0, \\
\sigma_H(y), & \text{if } \sigma_H(x) = 0,
\end{cases}
\end{gather*}
where $\sigma(x)$ and $\sigma(y)$ are the sign sequences of $x$ and
$y$. This product has a geometric interpretation: $xy$ is the face
entered by moving a small positive distance along a straight line from
a point in $x$ towards a point in $y$. In the special case where $y$
is a chamber, the product $xy$ is the chamber that has $x$ as a face
and that is separated from $y$ by the fewest number of hyperplanes in
$\Arrangement$ \cite[{\S}2C]{BrownDiaconis1998}. It is straightfoward
to verify that this product gives $\Faces$ the structure of an
associative semigroup with identity, and that $x^2 = x$ and $xyx = xy$
for all $x,y \in \Faces$. Semigroups satisfying these identities are
called \defn{left regular bands}.

The semigroup algebra $k\Faces$ is called the \defn{face semigroup
algebra} of $\Arrangement$ over the field $k$. It consists of finite
$k$-linear combinations of elements of $\Faces$ with multiplication
extended $k$-linearly from the product of $\Faces$. 

The semigroup $\Faces$ is also a partially ordered set with respect to
the relation $x\leq y$ if and only if $xy=y$. Equivalently, $x \leq y$
if and only it $x \subseteq \overline{y}$, where $\overline{y}$
denotes the closure of the set $y$. Note that the chambers of the
arrangement are precisely the faces that are maximal with respect to
this partial order. If $x\leq y$, then we say that $x$ is a face of
$y$ or that $y$ contains $x$ as a face.

\subsubsection{Support map and intersection lattice}
\label{sss:IntersectionLattice}
For each face $x \in \Faces$, the \defn{support} of $x$, denoted by
$\supp(x)$, is the intersection of all hyperplanes in $\Arrangement$
that contain $x$. Equivalently, $\supp(x)$ is the subspace of $V$
spanned by the vectors in $x$. The \defn{dimension} of $x$ is the
dimension of the subspace $\supp(x)$. 

The \defn{intersection lattice} $\IntersectionLattice$ of
$\Arrangement$ is the image of $\supp$; that is, $\IntersectionLattice
= \supp(\Faces)$. The elements of $\IntersectionLattice$ are subspaces
of $V$ and are ordered by inclusion. (N.B. Some authors order
$\IntersectionLattice$ by reverse inclusion rather than inclusion.)
With this partial order, $\IntersectionLattice$ is a finite lattice,
where the meet ($\vee$) of two subspaces is their intersection, and
the join ($\wedge$) of two subspaces is the smallest subspace that
contains both. 

It is straightforward to show that $\supp(x) \leq \supp(y)$ for all
$x,y\in\Faces$ with $x \leq y$. Therefore, $\supp$ is an
order-preserving poset surjection. Moreover, $\supp(xy) = \supp(x)
\vee \supp(y)$ for all $x,y\in\Faces$, so $\supp$ is also a semigroup
homomorphism, where $\IntersectionLattice$ is viewed as a semigroup
with product $\vee$. The elements of $\Faces$ also satisfy $xy=x$ if
$\supp(x)\geq\supp(y)$. Proofs of these statements can be found in
\cite[Appendix A]{Brown2000}.

\subsubsection{Invariant subalgebra}
\label{sss:InvariantSubalgebra}

Since $W$ is a group of orthogonal transformations of the vector space
$V$, there is a natural action of $W$ on $V$: the action of $w \in W$
on $\vec v \in V$ is the image of $\vec v$ under the transformation
$w$. This action permutes the set $\Arrangement$
\cite[Proposition~1.2]{Humphreys1990}, so it induces an action of $W$
on $\IntersectionLattice$ and on $\Faces$. This induced action
preserves the semigroup structure of $\Faces$ and
$\IntersectionLattice$, so it extends linearly to an action on
$k\Faces$ and $k\IntersectionLattice$. 

Let $\InvariantSubalgebra$ denote the subalgebra of $k\Faces$
consisting of the elements of $k\Faces$ fixed by all elements of $W$: 
\begin{gather*}
\InvariantSubalgebra = \Big\{ a \in k\Faces : w(a) = a \text{ for all
} w\in W \Big\}.
\end{gather*}
The following was first proved by T. P. Bidigare \cite{Bidigare1997}.
Another proof was given by K. S. Brown \cite[Theorem 7]{Brown2000}.

\begin{Theorem}[T. P. Bidigare]
\label{t:Bidigare}
Let $W$ be a finite reflection group and let $k\Faces$ denote the face
semigroup algebra of the reflection arrangement of $W$. The
$W$-invariant subalgebra $\InvariantSubalgebra$ is anti-isomorphic to
the descent algebra $\DescentAlgebra$ of $W$.
\end{Theorem}

We briefly describe an anti-isomorphism. The faces of the fundamental
chamber $c$ are parametrized by the subsets of $S$: if $J \subseteq
S$, then there is a unique face $c_J$ of $c$ that is fixed by all
elements in $J$ \cite[{\S}I.5F]{Brown1989}. Furthermore, every face of
$\Arrangement$ is in the $W$-orbit of a unique face of $c$
\cite[{\S}I.5F]{Brown1989}. So if $\Orbit_J$ denotes the $W$-orbit of
$c_J$, then the elements $\bs x_J = \sum_{y \in \Orbit_J} y$ form a
basis of $\InvariantSubalgebra$. The map defined by sending $\bs
x_J$ to $x_J$ is an anti-isomorphism from $\InvariantSubalgebra$ onto
$\DescentAlgebra$.

\subsection{The Symmetric Group} 
\label{ss:SymmetricGroupExample}
We describe the above ideas in combinatorial terms for the symmetric
group $\Sn$. The results in this section are not crucial to what
follows, and will only be used in the proof of Theorem
\ref{t:QuiverOfDesAlgA} to give a combinatorial description of the
quiver of the descent algebra $\DescentAlgebra[\Sn]$.

For $n\in\mathbb N$, let $[n] = \{1, \ldots, n\}$. A \defn{set
partition} of $[n]$ is a collection of nonempty subsets $B = \{B_1,
\ldots, B_r\}$ of $[n]$ such that $\bigcup_i B_i = [n]$ and $B_i \cap
B_j = \emptyset$ for $i \neq j$. The sets $B_i$ in $B$ are called the
\defn{blocks} of $B$. A \defn{set composition} of $[n]$ is an ordered
set partition of $[n]$, which we denote by $(B_1, \ldots, B_r)$. An
\defn{integer partition} of $n\in \mathbb N$ is a collection of
positive integers that sum to $n$.

\subsubsection{Braid arrangement}
Fix $n \in \mathbb N$. The \defn{braid arrangement} is the hyperplane
arrangement $\Arrangement$ in $V = \mathbb R^n$ consisting of the
hyperplanes $H_{ij} = \{ \vec v \in V : v_i = v_j \}$ for $1 \leq i <
j \leq n$. The group of transformations generated by the reflections
in the hyperplanes in $\Arrangement$ is identified with the symmetric
group acting on $V$ by permuting coordinates: $\omega( v_1, \ldots,
v_n ) = \left(v_{\omega^{-1}(1)}, \ldots, v_{\omega^{-1}(n)}\right)$
for $\omega \in \Sn$ and $v \in V$. The reflections in the hyperplanes
$H_{i,j}$ correspond to the transpositions $(i,j)\in\Sn$.

\subsubsection{Faces}
Let $\vec v = (v_1,v_2,\ldots,v_n) \in \mathbb R^n$ be a vector in a
chamber $c$ of the braid arrangement $\Arrangement$. Then $\vec v$ is
not on any of the hyperplanes $H_{ij}$, so all the coordinates of
$\vec v$ are distinct. Therefore, there exists $\omega \in
\Sn$ such
that $v_{\omega(1)} < \cdots < v_{\omega(n)}$. All vectors in $c$
satisfy this identity, so $c$ can be identified with the permutation
$\omega$. The faces of $c$ are obtained by changing some of the
inequalities to equalities, so the \defn{faces} of $\Arrangement$ can
be identified with set compositions of $[n]$. For example, the set
composition $( \{5\}, \{1,3,4\}, \{2,6\} )$ is identified with the
face $\{\vec v \in \mathbb R^5 : v_5 < v_1 = v_3 = v_4 < v_2=v_6 \}
=H_{1,5}^-\cap H_{1,3} \cap H_{3,4} \cap H_{2,4}^- \cap H_{2,6}$,
where $H_{i,j}^-=\{\vec v: v_i>v_j\}$.

The \defn{partial order} is given by $(B_1, \ldots, B_m) \leq (C_1,
\ldots, C_l)$ if and only if $(C_1, \ldots, C_l)$ consists of a set
composition of $B_1$, followed by a set composition of $B_2$, and so
forth. The \defn{action} of $\Sn$ on set compositions is given by
permuting the underlying set: $\omega(B_1, \ldots, B_r) = (
\omega(B_1), \ldots, \omega(B_r) )$. And if $(B_1, \ldots, B_l)$ and
$(C_1, \ldots, C_m)$ are set compositions of $[n]$, their
\defn{product} is the set composition of $[n]$ given by 
the formula
\begin{align*}
(B_1, ..., B_l) & (C_1, ..., C_m) \\
& = (B_1 \cap C_1, ..., B_l \cap C_1, 
\ \ldots, \
B_1 \cap C_m, ..., B_l \cap C_m)^{\scissors},
\end{align*}
where $\scissors$ means ``delete empty intersections''. 

\subsubsection{Intersection lattice}
The elements of the intersection lattice $\IntersectionLattice$ of
$\Arrangement$ are identified with set partitions of $[n]$ via the
following bijection,
\begin{align*}
\{B_1, & \ldots, B_r\} \leftrightarrow \\[-2ex]
&\Big\{\vec v \in V : v_i = v_j\text{ if } i,j \in B_h \text{ for some
}h\in[r] \Big\}
= \bigcap_{h=1}^{r} \left(\bigcap_{i,j \in B_h} H_{ij}\right),
\end{align*}
where $\{B_1, \ldots, B_r\}$ is a set partition of $[n]$. 

Under this identification, if $B$ and $C$ are set partitions of $[n]$,
then $B \lessdot C$ if and only if $B$ is obtained from $C$ by merging
two blocks of $C$. The action of $\Sn$ on $\IntersectionLattice$ is
given by $\omega(\{B_1, \ldots, B_r\}) = \{ \omega(B_1), \ldots,
\omega(B_r) \}$. The \defn{support map} sends a set composition
$(B_1,\ldots,B_m)$ to the underlying set partition
$\{B_1,\ldots,B_m\}$. 

The $\Sn$-orbit of a set partition $\{B_1,\ldots,B_m\}$ of $[n]$
depends only on the sizes of the blocks $B_i$, so
$\IntersectionLattice/\Sn$ can be identified with the poset of integer
partitions of $n$. Under this identification, for any two integer
partitions $p$ and $q$ of $n$, we have $p \lessdot q$ if and only if
$p$ is obtained from $q$ by adding two elements of $q$.

\section{The Quiver of a Split Basic Algebra}
\label{s:QuiverOfSplitBasicAlgebra}

This section recalls definitions and results from the theory of finite
dimensional algebras. Our main references are
\cite{AuslanderReitenSmalo1995,Benson1998:I,Assem2006}.

Let $k$ be a field and $A$ a finite dimensional $k$-algebra. 
An element $a \in A$ is an \defn{idempotent} if $e^2 = e$. Two
idempotents $e,f\in A$ are \defn{orthogonal} if $ef = 0 = fe$. An
idempotent $e\in A$ is \defn{primitive} if it cannot be written as $e
= f + g$ with $f$ and $g$ non-zero orthogonal idempotents of $A$. A
\defn{\CompleteSystemPunctuated} of $A$  is a set $\{e_1, e_2, \ldots,
e_n\}$ of primitive idempotents of $A$ that are pairwise orthogonal
and that sum to $1\in A$.

The \defn{Jacobson radical} of $A$ is the smallest ideal $\rad(A)$ of
$A$ such that $A/\mathord{\rad(A)}$ is semisimple. If
$A/\mathord{\rad(A)}$ is isomorphic, as a $k$-algebra, to a direct
product of copies of $k$, then $A$ is said to be a 
\defn{split basic algebra}. Equivalently, $A$ is a split basic algebra
if and only if all the simple $A$-modules are one dimensional.

The \defn{quiver} of a split basic $k$-algebra $A$ is the directed
graph $Q$ constructed as follows. Let $\{e_v: v \in \mc V\}$ be a
\CompleteSystem of $A$, where $\mc V$ is some index set. There is one
vertex $v$ in $Q$ for each idempotent $e_v$ in $\{e_v: v \in \mc V\}$.
If $x,y \in \mc V$, then the number of arrows in $Q$ from $x$ to $y$
is $\dim_k e_y\big(\rad(A)/\mathord{\rad^2(A)}\big)e_x$. This
construction does not depend on the \CompleteSystem 
(see \cite[Definition~4.1.6]{Benson1998:I} or \cite[Lemma~II.3.2]{Assem2006}).

If $\alpha$ is an arrow in a quiver (directed graph) beginning at a
vertex $x$ and ending at a vertex $y$, then we write
$x\labelledarrow{\alpha}y$. If there is exactly one arrow from $x$ to
$y$, then we drop the label and write $x\arrow y$.
The \defn{path algebra} $kQ$ of a quiver $Q$ is the $k$-algebra with
basis the set of paths in $Q$ and with multiplication defined on paths
by 
\begin{gather*}
(w_0 \labelledarrow{\alpha_1}\cdots\labelledarrow{\alpha_s} w_s) \cdot
(v_0 \labelledarrow{\beta_1}\cdots\labelledarrow{\beta_2} v_r) = 
\begin{cases}
(v_0 \labelledarrow{\beta_1}\cdots\labelledarrow{\beta_2} v_r 
\labelledarrow{\alpha_1} w_1\labelledarrow{\alpha_2}\cdots\labelledarrow{\alpha_s} w_s),
& \text{if }w_0 = v_r, \\
\hfill0\hfill,
& \text{if }w_0 \neq v_r,
\end{cases}
\end{gather*}
where 
$(w_0 \labelledarrow{\alpha_1}\cdots\labelledarrow{\alpha_s} w_s)$
and 
$(v_0 \labelledarrow{\beta_1}\cdots\labelledarrow{\beta_2} v_r)$
are paths in $Q$. 

Let $F$ denote the ideal in $kQ$ generated by the arrows of $Q$. An
ideal $I \subseteq kQ$ is said to be \defn{admissible} if there exists
an integer $m\geq2$ such that $F^m \subseteq I \subseteq F^2$. This
notion is useful for identifying the quiver of a split basic
$k$-algebra as the following result demonstrates
\cite[Theorem~III.1.9(d)]{AuslanderReitenSmalo1995}.

\begin{Theorem}
\label{t:ARS}
$Q$ is the quiver of a finite dimensional split basic $k$-algebra $A$
if and only if $A \cong kQ/I$, where $I$ is an admissible ideal of
$kQ$. In particular, if $\varphi:kQ\twoheadrightarrow A$ is a
surjection of $k$-algebras with an admissible kernel, then $Q$ is the
quiver of $A$.
\end{Theorem}

The following result will be helpful to define $k$-algebra morphisms.
A proof can be found in \cite[Theorem II.1.8]{Assem2006}.
\begin{Theorem}
\label{t:ExtendingMapsToMorphisms}
Let $Q$ be a finite quiver and $A$ a finite dimensional $k$-algebra.
If $f$ is a function from the set of vertices and arrows of $Q$ into
$A$ such that
\begin{enumerate}
\item $\sum_v f(v)=1$, $f(v)^2=f(v)$ and $f(u)f(v)=0$ for all
 vertices $u,v$, and
\item $f(u\arrow v)=f(v)f(u\arrow v)f(u)$ for every arrow $u\arrow v$,
\end{enumerate}
then there exists a unique $k$-algebra homomorphism
$\varphi:kQ\to A$ such that $\varphi(v)=f(v)$ and
$\varphi(u\arrow v)=f(u\arrow v)$ for all vertices $v$ and
all arrows $u\arrow v$ of $Q$.
\end{Theorem}

\section{$k\Faces$ and $\InvariantSubalgebra$ are Split Basic Algebras}
\label{s:MyAlgebrasAreBasic}

This section establishes that $k\Faces$ and $\InvariantSubalgebra$ are
split basic algebras. That $\InvariantSubalgebra$ is a split basic
algebra follows from various sources since the irreducible
representations of the descent algebra are known to be one dimensional
(see, for example, \cite[Theorem 3]{Solomon1976}). We give a proof
based on \cite{Bidigare1997,Brown2000}. 

\begin{Proposition}
\label{p:MyAlgebrasAreBasic}
$k\Faces$ and $\InvariantSubalgebra$ are split basic algebras.
\end{Proposition}
\begin{proof}
We begin by showing that $k\Faces$ is a split basic algebra. As
mentioned in \S\ref{sss:IntersectionLattice}, the support map $\supp: \Faces
\to \IntersectionLattice$ is a surjective semigroup homomorphism.
Therefore, it extends linearly to a surjective $k$-algebra
homomorphism $\supp: k\Faces \to k\IntersectionLattice$. The algebra
$k\IntersectionLattice$ is isomorphic to the $k$-algebra
$\prod_{X\in\IntersectionLattice} k$. Indeed, the elements defined
recursively by the formula $E_X=X-\sum_{Y>X}E_Y$, one for each
$X\in\IntersectionLattice$, form a basis and a \CompleteSystem for
$k\IntersectionLattice$ \cite{Solomon1967}. Since the kernel of
$\supp$ is nilpotent, standard ring theory implies that
$\ker(\supp)=\rad(k\Faces)$. It follows that $k\Faces$ is a split
basic algebra. 

Since $\supp: k\Faces \twoheadrightarrow k\IntersectionLattice$ is a
surjective $W$-equivariant algebra homomorphism, it restricts to an
algebra surjection $\InvariantSubalgebra \twoheadrightarrow
(k\IntersectionLattice)^W$, where $(k\IntersectionLattice)^W$ is the
$W$-invariant subalgebra of $k\IntersectionLattice$. Let $E_X$ be the
elements defined above. Since $w(E_X)=E_{w(X)}$ for all $w\in W$ and
$X\in\IntersectionLattice$, it follows that the elements
$\sum_{X\in\Orbit} E_X$, one for each $W$-orbit $\Orbit$ of elements
of $\IntersectionLattice$, form a basis and a \CompleteSystem for
$(k\IntersectionLattice)^W$. Thus,
$(k\IntersectionLattice)^W\cong\prod_{\Orbit\in\IntersectionLattice/W}k$.
Since $(k\IntersectionLattice)^W$ is semisimple and
$\ker(\supp|_{\InvariantSubalgebra})$ is nilpotent (because
$\ker(\supp)$ is), it follows that the radical of
$\InvariantSubalgebra$ is $\ker(\supp|_{\InvariantSubalgebra})$. Thus,
$\InvariantSubalgebra$ is a split basic algebra.
\end{proof}

\section{Complete Systems of Primitive Orthogonal Idempotents}
\label{s:CompleteSystems}

In this section we construct a \CompleteSystem for $k\Faces$
that is permuted by the elements of $W$. This allows us to construct a
\CompleteSystem for $\InvariantSubalgebra$. 
A \CompleteSystem for $\DescentAlgebra$ was constructed previously
\cite{BergeronBergeronHowlettTaylor1992}, but the construction
presented here is new and better suited to our needs because of the
close relationship between the two systems.

For each $X \in \IntersectionLattice$ let $\Orbit_X = \{ w(X) : w \in
W \}$ denote the $W$-orbit of $X$. These orbits form a poset
$\IntersectionLattice/W = \{ \Orbit_X : X \in \IntersectionLattice \}$
with partial order given by $\Orbit_X \leq \Orbit_Y$ if and only if
there exists $w \in W$ with $w(X) \leq Y$. 

\begin{Remark}
\label{r:IntLatAndSubsetsOfS}
The poset $\IntersectionLattice/W$ is isomorphic to a poset of
equivalence classes of subsets of $S$. Indeed, define a relation on
subsets $J,K \subseteq S$ by setting $J \sim K$ if and only if
$\supp(c_J)$ and $\supp(c_K)$ belong to the same $W$-orbit, where
$c_J$ and $c_K$ are the largest faces of the fundamental chamber
$c$ that are fixed by $J$ and $K$, respectively.
Equivalently, $J \sim K$ if and only if $W_J$ and $W_K$ are conjugate
subgroups of $W$. The poset $S/\mathord{\sim}$, with partial order
induced by reverse inclusion of subsets of $S$, is isomorphic to
$\IntersectionLattice/W$. 
\end{Remark}

\begin{Theorem}
\label{t:CompSystem}
For each $X\in\IntersectionLattice$, fix a linear combination $\ell_X$
of faces of support $X$ whose coefficients sum to $1$ and suppose
that they satisfy the identity
\begin{align}
\label{e:InvariantElls}
w(\ell_X) = \ell_{w(X)} \text{ for all }
w\in W, X\in\IntersectionLattice.
\end{align}
Then the elements defined recursively using the equation
\begin{align}
\label{e:IdempotentsFormula}
e_X = \ell_X - \sum_{Y>X} \ell_X e_Y
\end{align}
one for each $X\in\IntersectionLattice$, form a \CompleteSystem for
$k\Faces$, and they satisfy $w(e_X)=e_{w(X)}$ for every $w\in W$ and
$X\in\IntersectionLattice$. The elements
\begin{align}
\label{e:DesAlgIdemps}
\varepsilon_\Orbit = \sum_{X\in\Orbit} e_X,
\end{align}
one for each $\Orbit\in\IntersectionLattice/W$, form a \CompleteSystem
for $\InvariantSubalgebra$.
\end{Theorem}

Examples of elements $\ell_X$ satisfying the above hypotheses will be
presented below.

\begin{proof}
In \cite[Theorem 4.2]{Saliola2007a} and \cite[Theorem
5.2]{Saliola2008a} it was shown that the elements $e_X$ form a
\CompleteSystem for $k\Faces$. (This is proved by first establishing
Lemma \ref{l:IdempotentLemma} below and inducting on the codimension
of $X$ in $V$.) Induction on the codimension of
$X\in\IntersectionLattice$ establishes that $w(e_X) = e_{w(X)}$ for
all $w\in W$ and all $X \in \IntersectionLattice$. Therefore, the
elements $\sum_{Y \in \Orbit} e_Y$ are invariant under the action of
$W$, so they belong to $\InvariantSubalgebra$. They are orthogonal
idempotents since sums of orthogonal idempotents are again orthogonal
idempotents. They sum to $1$ since $\sum_{X\in\IntersectionLattice}
e_X = 1$. Finally, they are primitive because there are enough of
them: the number of elements in a \CompleteSystem for a split basic
algebra $A$ is the dimension of $A/\rad(A)$ (this follows from
\cite[Corollary 1.7.4]{Benson1998:I}), which in this case is
$|\IntersectionLattice/W|$ by the proof of Proposition
\ref{p:MyAlgebrasAreBasic}. 
\end{proof}

The idempotents $e_X$ satisfy the following remarkable property that
we will use on occasion. A proof can be found in 
\cite[Lemma 4.1]{Saliola2007a} and \cite[Lemma 5.1]{Saliola2008a}.
\begin{Lemma}[\cite{Saliola2007a,Saliola2008a}]
\label{l:IdempotentLemma}
Let $y \in \Faces$ and $X \in \IntersectionLattice$. If $\supp(y)
\not\leq X$, then $y e_X = 0$.
\end{Lemma}

Next we present some examples of elements $\ell_X$ satisfying the
above hypotheses. 

\subsection{First Complete System}
\label{ss:FirstSystem}
For each $X \in \IntersectionLattice$, let $\ell_X$ denote the
normalized sum of all faces of support $X$:
\begin{gather*}
\ell_X = \frac{1}{{}^\#\{x \in \Faces : \supp(x) = X\}}
 \left(\sum_{\supp(x) = X} x\right).
\end{gather*}
Then $w(\ell_X) = \ell_{w(X)}$ for all $w \in W$ and
$X\in\IntersectionLattice$.

\subsection{Second Complete System}
\label{ss:SecondSystem}

For every orbit $\Orbit\in\IntersectionLattice/W$, fix a face
$f_\Orbit\in\Faces$ such that $\supp(f_\Orbit)\in\Orbit$. For each
$X\in\IntersectionLattice$, let $f_X = f_{\Orbit_X}$ and define
\begin{gather}
\label{e:SecondSystem}
\ell_X = 
\frac1{L_X}
 \left( \sum_{ z \in \Orbit_{f_X} \atop \supp(z) = X } z \right),
\quad \text{ where }
L_X 
=\Big|\{ z \in \Orbit_{f_X} : \supp(z) = X \}\Big|.
\end{gather}
Note that $L_X$ is the index of the stabilizer subgroup $W_x$ of $x$,
where $x$ is any face of support $X$, in the stabilizer subgroup $W_X$
of $X$. It follows that every $w\in W$ induces a bijection between 
\begin{gather*}
T_X = \{z\in\Orbit_{f_X} : \supp(z)=X\}
\,\text{ and }\,
T_{w(X)} = \{z\in\Orbit_{{f_{w(X)}}} : \supp(z)=w(X)\},
\end{gather*}
so $w(\ell_X)=\ell_{w(X)}$ for all $X\in\IntersectionLattice$ and all
$w\in W$.

\subsection{Third Complete System}
\label{ss:ThirdSystem}

If $(W,S)$ is a Coxeter system with fundamental chamber $c$, then the
faces of $c$ are parametrized by the subsets of $S$: if $J \subseteq
S$, then there is a unique largest face $c_J$ of the fundamental
chamber $c$ that is fixed by all elements of $J$
\cite[{\S}I.5F]{Brown1989}. For $J\subseteq S$, let $\bs x_J$ denote
the sum of the faces in the $W$-orbit of $c_J$ (see also
\S\ref{sss:InvariantSubalgebra}). For each orbit
$\Orbit\in\IntersectionLattice/W$, fix a subset $J_\Orbit\subseteq S$
such that $\supp(c_{J_\Orbit}) \in \Orbit$ and define numbers
\begin{align*}
L_\Orbit=
\big|\{z\in\Orbit_{x_{J_\Orbit}}:\supp(z)=\supp(x_{J_\Orbit})\}\big|.
\end{align*}

\begin{Proposition}
\label{p:ThirdSystem}
The elements $\varepsilon_\Orbit$, one for each
$\Orbit\in\IntersectionLattice/W$, defined recursively by the formula
\begin{gather*}
\varepsilon_\Orbit 
  = \frac1{L_\Orbit} \bs x_{J_\Orbit} 
  - \sum_{\Orbit'>\Orbit} 
   \left(\frac1{L_\Orbit} \bs x_{J_\Orbit} \right) 
   \varepsilon_{\Orbit'},
\end{gather*}
form a \CompleteSystem for $\InvariantSubalgebra$.
\end{Proposition}
\begin{proof}
For each $\Orbit\in\IntersectionLattice/W$, let $f_\Orbit =
x_{J_\Orbit}$, and let $f_X = f_{\Orbit_X}$ for each
$X\in\IntersectionLattice$. Define $\ell_X$ using Equation
\eqref{e:SecondSystem}. Then an induction on the corank of $\Orbit$
establishes that the elements defined by Equation \eqref{e:DesAlgIdemps}
are equal to the elements $\varepsilon_\Orbit$ defined above.
\end{proof}

\begin{Remark}
\label{r:IdempotentsInDescentAlgebra}
Proposition \ref{p:ThirdSystem} leads to a construction of a
\CompleteSystem directly within the descent algebra \DescentAlgebra.
Let $S/\mathord{\sim}$ denote the poset defined in Remark
\ref{r:IntLatAndSubsetsOfS}. For each $\Orbit \in S/\mathord{\sim}$,
fix a subset $J_\Orbit\subseteq S$ with $J_\Orbit \in \Orbit$ and
define elements $\varepsilon_\Orbit$, one for each $\Orbit\in
S/\mathord{\sim}$, recursively by the formula
\begin{gather*}
\varepsilon_\Orbit 
  = \frac1{L_\Orbit} x_{J_\Orbit} 
  - \sum_{\Orbit'>\Orbit} \varepsilon_{\Orbit'}
     \left(\frac1{L_\Orbit} 
     x_{J_\Orbit} \right),
\end{gather*}
where $x_{J_\Orbit}$ are the basis elements of $\DescentAlgebra$
as defined in \S\ref{ss:DescentAlgebra} and where
$L_\Orbit$ is the index of $W_J$ in the normalizer of $W_J$.
\end{Remark}

\begin{Remark}
The construction of the idempotents $\varepsilon_\Orbit$ in Remark
\ref{r:IdempotentsInDescentAlgebra} is very similar to the
construction of the idempotents $e_X$ in Theorem \ref{t:CompSystem}.
In both cases we start with an algebra $A$ and a function $s$ from a
basis of $A$ to a poset $P$. Using $s$ and $P$, a \CompleteSystem is
constructed for $A$ using Equation \eqref{e:IdempotentsFormula}. It
would be interesting to determine conditions on $A$, $P$, and $s$ to
ensure that this construction provides a \CompleteSystem for $A$.
\end{Remark}

\section{A $W$-Equivariant Surjection}
\label{s:WEquivariantSurjection}

In this section we define a quiver $\Quiver$ and a $W$-equivariant
surjection $\varphi:k\Quiver\to k\Faces$ of $k$-algebras. We use this
homomorphism in later sections to deduce properties of the quiver of
$\InvariantSubalgebra$. Recall that we write $Y\lessdot X$ if and only
if $Y<X$ and there exists no $Z\in\IntersectionLattice$ such that
$Y<Z<X$.

\begin{Definition}
Let $\Quiver$ be the directed graph on the vertex set
$\IntersectionLattice$ and with exactly one arrow $X \arrow Y$ if and
only if $Y \lessdot X$. 
\end{Definition}

In \cite[Corollary 8.4]{Saliola2008a} it is shown that $\Quiver$ is
the quiver of $k\Faces$, which will also follow from the theorem below
(see Corollary \ref{c:QuiverOfSemigroupAlgebra}). In that article it
is also shown that $k\Faces$ is a Koszul algebra, but this fact won't
be necessary here. 

\begin{Theorem}
\label{t:WEquivariantPhi}
Let $\{e_X\}_{X\in\IntersectionLattice}$ denote a \CompleteSystem for
$k\Faces$ as defined in Theorem \ref{t:CompSystem}. 
Fix an orientation $\epsilon_X$ on each subspace
$X\in\IntersectionLattice$ and define numbers $[x:y]$ for pairs of faces
satisfying $x\lessdot y$ by
\begin{align}
\label{e:IncidenceNumbers}
[x:y] &= \epsilon_{\supp(x)}(\vec x_1, \ldots, \vec x_d)
 \epsilon_{\supp(y)}(\vec x_1, \ldots, \vec x_d, \vec y_1),
\end{align}
where $\vec x_1, \ldots, \vec x_d$ is a basis of $\supp(x)$ and $\vec
y_1$ is a vector in $y$.

Let $\varphi$ be the function defined on the vertices and arrows of
$\Quiver$ by 
\begin{align*}
\varphi(X)=e_X \qquad\text{and}\qquad
 \varphi(X\arrow Y) = \ell_Y\Big([y:x] x + [y:x']x'\Big)e_X,
\end{align*}
where $y$ is any face of support $Y$ and where $x$ and $x'$ are the two
faces of support $X$ having $y$ as a face. Then $\varphi$ extends
uniquely to a surjection of $k$-algebras
$\varphi:k\Quiver\twoheadrightarrow k\Faces$, the kernel of $\varphi$
is generated as an ideal by the sum of all the paths of length two in
$\Quiver$, and $\varphi$ is $W$-equivariant with respect to the
following action of $W$ on $k\Quiver$:
\begin{align*}
w \Path Xt
= \sigma_{X_0}(w)\sigma_{X_t}(w)\Big(w(X_0) \arrows w(X_t)\Big),
\end{align*}
where $\sigma_X(w)$, for $X\in\IntersectionLattice$ and $w\in
W$, is defined by the equation
\begin{align}
\label{e:OrientationNumbers}
\sigma_X(w) &= \epsilon_X\Big(\vec x_1, \ldots, \vec x_d\Big)
 \epsilon_{w(X)}\Big(w(\vec x_1), \ldots, w(\vec x_d)\Big),
\end{align}
where $\vec x_1, \ldots, \vec x_d$ is a basis of $X$.
\end{Theorem}

We will prove this by a sequence of lemmas. But before we do, let us
record a few properties of the numbers defined in Equations
\eqref{e:IncidenceNumbers} and \eqref{e:OrientationNumbers}.

It is straightforward to prove that the \defn{incidence numbers}
defined in Equation \eqref{e:IncidenceNumbers} satisfy the identity 
\begin{align}
\label{e:IncidenceNumberIdentityI}
[x:y]=[x':x'y], \text{ if } x,x'\in\Faces \text{ and } \supp(x')=\supp(x).
\end{align}
They also satisfy the following identity,
\begin{align}
\label{e:IncidenceNumberIdentityII}
[z:y][y:x]+[z:u][u:x]=0,
\end{align}
where $y$ and $u$ are the two faces in the interval
$\{f\in\Faces:z\lessdot f\lessdot x\}$. A proof of this can be found
in \cite[Lemma 2 in \S5C]{BrownDiaconis1998}, 

\begin{Remark}
The incidence numbers were defined by Kenneth S. Brown and Persi
Diaconis who used them to compute the multiplicities of the
eigenvalues of random walks on the chambers of a hyperplane
arrangement. The numbers get their name from the fact that they form a
system of ``incidence numbers'' in the sense of homology theory of
regular cell complexes. See \cite[\S5]{BrownDiaconis1998} for details.
\end{Remark}

The number $\sigma_X(w)$ defined in Equation
\eqref{e:OrientationNumbers} measures whether $w$ maps a positively
oriented basis of $X$ to a positively or negatively oriented basis of
$w(X)$. Note that if $w(X) = X$, then $\sigma_X(w)$ is $1$ if and only
if the restriction of $w$ to $X$ is orientation-preserving, and is
$-1$ otherwise. In particular, if $w(X)=X$, then the number
$\sigma_X(w)$ does not depend on the choice of $\epsilon_X$. And since
$w$ is an orthogonal transformation of $V$: if $w(X) = X$, $w(Y) = Y$
and $Y\lessdot X$, then $\sigma_X(w) \sigma_Y(w) = -1$ if and only if
$w$ interchanges the two halfspaces of $X$ determined by $Y$.

\subsection{Proof of Theorem \ref{t:WEquivariantPhi}}

We begin by showing that $\varphi$ is well-defined.

\begin{Lemma}
$\varphi:k\Quiver\to k\Faces$ is a well-defined homomorphism of
$k$-algebras.
\end{Lemma}

\begin{proof}
There are a three issues that need to be addressed with the definition
of $\varphi$. First is the fact that there are exactly two faces of
support $X$ having $y$ as a face, which is a well known result
\cite[\S I.4E Proposition 3]{Brown1989}. The second issue is the claim
that $\varphi(X\arrow Y)$ does not depend on the choice of $y$.
Indeed, since $\ell_Y$ is a linear combination of faces of support
$Y$, we have $\ell_Yy'=\ell_Y$ for any face $y'$ with $\supp(y')=Y$
(this is because $yy'=y$ if $\supp(y)\geq\supp(y')$; see
\S\ref{sss:IntersectionLattice}). So if $y'$ is another face with $\supp(y')=Y$,
then
\begin{align*}
\ell_Y\Big([y:x] x + [y:x']x'\Big)e_X
&= \ell_Y\Big([y:x] y'x + [y:x']y'x'\Big)e_X \\
&= \ell_Y\Big([y':y'x] y'x + [y':y'x']y'x'\Big)e_X,
\end{align*}
where we used Equation \eqref{e:IncidenceNumberIdentityI} to obtain
the last equality. Since $y'x$ and $y'x'$ are the two faces of support
$X$ having $y'$ as a face, the claim follows.

The third issue is that $\varphi$ extends uniquely to a homomorphism
of $k$-algebras. Since the images of the vertices form a
\CompleteSystemPunctuated, it suffices to show that
$\varphi(Y)\varphi(X\arrow Y)\varphi(X)=\varphi(X\arrow Y)$ for all
arrows $X\arrow Y$ in $\Quiver$ (Theorem
\ref{t:ExtendingMapsToMorphisms}). Since $\varphi(X)=e_X$ is an
idempotent for each vertex $X\in\IntersectionLattice$, it follows
immediately that $\varphi(X\arrow Y)\varphi(X)=\varphi(X\arrow Y)$. It
remains to show that $\varphi(Y)\varphi(X\arrow Y)=\varphi(X\arrow
Y)$. Using Equation \eqref{e:IdempotentsFormula} we write,
\begin{align*}
\varphi(Y)\varphi(X\arrow Y) 
&= \left(\ell_Y - \sum_{U>Y}\ell_Y e_U\right) \ell_Y
   \Big([y:x] x + [y:x']x'\Big) e_X \\
&= \varphi(X\arrow Y) - 
   \sum_{U>Y}\ell_Y e_U \Big([y:x] x + [y:x']x'\Big) e_X.
\end{align*}
We will show this is $\varphi(X\arrow Y)$ by showing each
term in the summation is zero. 

Suppose $U>Y$. Since $Y\lessdot X$, either $U=X$ or $U\not\leq X$. If
$U=X$, then let $u$ be a face of support $U$ and note that $ux=u=ux'$
because $\supp(u)\geq\supp(x)$ (\S\ref{sss:IntersectionLattice}). Since
$e_U$ is a linear combination of elements of support at least $U=X$, we
have $e_Uu=e_U$. Thus,
\begin{align*}
e_U \Big([y:x] x + [y:x']x'\Big)=e_U\Big([y:x]u+[y:x']u\Big)=0
\end{align*}
since $[y:x]=-[y:x']$. If $U\not\leq X$, then the fact that
$(e_Ua)u=e_Ua$ for any $a\in k\Faces$ and any $u\in\Faces$ with
$\supp(u)=U$ implies
\begin{align*}
e_U \Big([y:x] x + [y:x']x'\Big) e_X = e_U\Big([y:x]x+[y:x']x'\Big)(ue_X)=0, 
\end{align*}
where the last equality follows from Lemma \ref{l:IdempotentLemma}.
\end{proof}

\begin{Lemma}
\label{l:PhiIsSurjective}
$\varphi:k\Quiver\to k\Faces$ is surjective.
\end{Lemma}

\begin{proof}
Since the elements $e_X$ are orthogonal idempotents, to show that
$\varphi$ is surjective it suffices to show that $k\Faces e_X$ is in
the image of $\varphi$ for all $X\in\IntersectionLattice$. It follows
from Lemma \ref{l:IdempotentLemma} and Equation
\eqref{e:IdempotentsFormula} that the following is a basis of $k\Faces
e_X$:
\begin{align*}
\{xe_{X}: x\in\Faces, \supp(x)=X\}.
\end{align*}
A proof of this can be found in \cite[Lemma 5.1]{Saliola2007a} and
\cite[Lemma 6.1]{Saliola2008a}.

We proceed by induction on the rank of $X$ in $\IntersectionLattice$.
If the rank is zero, then $X$ is the intersection of all the
hyperplanes in $\Arrangement$. There is only one face that has this
support, the identity element of $\Faces$. Thus, $k\Faces e_X
\subseteq \img(\varphi)$ since $e_X = \varphi(X)$. 

Suppose $k\Faces e_Y \subseteq \img(\varphi)$ for all $Y$ satisfying
$\rank(Y) < r$. Let $X\in\IntersectionLattice$ with $\rank(X)=r$. Let
$x$ and $x'$ be two faces of support $X$ that are separated by exactly
one subspace $Y$ of $X$ having codimension one. We will prove that
$(x-x')e_X\in\img(\varphi)$. 
Let $y$ denote the face of support $Y$ that is common to both $x$ and
$x'$. Then $x = yx$ and $x' = yx'$. Thus, up to a sign
$(x - x')e_X$ is equal to 
\begin{align*}
\pm (x - x')e_X 
&= \Big([y: x] x + [y: x'] x'\Big)e_X  \\
&= y \Big([y: x] x + [y: x'] x'\Big)e_X  \\
&= y\ell_Y \Big([y: x] x + [y: x'] x'\Big)e_X  \\
&= y \varphi(X\arrow Y)  \\
&= y \varphi(Y)\varphi(X\arrow Y)  \\
&= (ye_Y) \varphi(X\arrow Y).
\end{align*}
Here we used the identity $y\ell_Y = y$. Since $Y$ is a proper
subspace of $X$, $\rank(Y)<\rank(X)=r$. By the induction hypothesis,
$y e_Y\in\img(\varphi)$. Hence, $xe_X-x'e_X\in\img(\varphi)$ for every
pair $x,x'$ of faces of support $X$ sharing a common codimension
one face. 

For every pair of faces $x$ and $x'$ of support $x$, there exists a
sequence of faces $x_0=x, x_1, \ldots, x_d=x'$ of support $X$ such
that $x_{i-1}$ and $x_{i}$ share a common codimension one face for
each $1 \leq i \leq d$ \cite[Proposition 3 of \S{I.4E}]{Brown1989}, it
follows that $xe_X-x'e_X\in\img(\varphi)$ for any pair of faces $x,x'$
of support $X$. Since the sum of the coefficients of $\ell_X$ is
nonzero, the elements $xe_X-x'e_X$, where $x,x'\in\Faces$ and
$\supp(x)=\supp(x')=X$, together with $\ell_Xe_X=e_X$ span the
subspace $k\Faces e_X$. Since $e_X=\varphi(X)$, it follows that
$k\Faces e_X\subseteq\img(\varphi)$. Thus, $\varphi$ is surjective.
\end{proof}

\begin{Lemma}
\label{l:KernelOfPhi}
The kernel of $\varphi:k\Quiver\to k\Faces$ is generated as an ideal
by the sum of all the paths of length two in $\Quiver$.
\end{Lemma}

\begin{proof}
Let $\rho$ be the sum of all the paths of length two in $\Quiver$ and
let $\mc I$ denote the ideal generated by $\rho$. If $X$ and $Y$ are
two vertices of $\Quiver$ with $Y\leq X$ and $\dim(X)=\dim(Y)+2$, then
$Y\rho X$ is the sum of all the paths of length two that begin at $X$
and end at $Y$. We begin by showing that these elements are in
$\ker(\varphi)$.

Suppose $(X\arrow Y\arrow Z)$ is a path of length two in $\Quiver$.
Let $z$ be a face of support $Z$ and $y$ a face of support $Y$.
Since $\supp(zy)=\supp(z)\vee\supp(Y)$, it follows that $zy$ has
support $Y$. By replacing $y$ with $zy$, we can suppose that
$z\lessdot y$. Thus,
\begin{align*} 
 \varphi(Y \arrow Z) \varphi(X \arrow Y) 
 =\ell_Z\Big([z:y]y+[z:y']y'\Big)\ell_Y\Big([y:x]x+[y:x']x'\Big)e_X,
\end{align*}
where $y$ and $y'$ are the two faces of support $Y$ having $z$ as a
face, and $x$ and $x'$ are the two faces of support $X$ having $y$ as
a face. Since $y$ and $y'$ have support $Y$, $y\ell_Y=y$ and
$y'\ell_Y=y'$. Thus,
\begin{align} 
\label{e:2PathImage}
 \varphi(Y \arrow Z) \varphi(X \arrow Y) 
 =\ell_Z\Big([z:y]y+[z:y']y'\Big)\Big([y:x]x+[y:x']x'\Big)e_X,
\end{align}
So $\varphi(X\arrow Y\arrow Z)$ is a linear combination of elements of
the form $\tilde x e_X$ with $\tilde x$ having support $X$, having a
codimension one face of support $Y$, and having a codimension two face
occuring in $\ell_Z$ with a nonzero coefficient. 

Let $z$ be a face occuring in $\ell_Z$ with a nonzero coefficient.
There are exactly two codimension one faces of $\tilde x$ that contain
$z$ as a face 
\cite[Lemma 2 of \S5C]{BrownDiaconis1998}; call them $y$ and $u$. So
$\tilde x$ can only appear in $\varphi(X\arrow Y \arrow Z)$ and
$\varphi(X\arrow U \arrow Z)$, where $Y=\supp( y)$ and $U=\supp(u)$.
Moreover, in Equation \eqref{e:2PathImage} exactly one of $y x$ or $y
x'$ can be $\tilde x$; we can suppose that $yx=\tilde x$. So, $\tilde
x e_X$ appears in $\varphi(X \arrow Y \arrow Z)$ with coefficient
$[z:y][y:\tilde x]$. Similarly, $\tilde
xe_X$ appears in $\varphi(X \arrow U \arrow Z)$ with coefficient
$[z:u][u:\tilde x]$. It follows from Equation
\eqref{e:IncidenceNumberIdentityII} that the coefficient of $\tilde
xe_X$ in the sum $\sum \varphi(X \arrow Y \arrow Z)$ is
$[z:y][y:\tilde x]+[z:u][u:\tilde x]=0$.

The above shows that $\mc I\subseteq\ker(\varphi)$. Let $X$ and $Y$ be
two two vertices in $\Quiver$, and let $M_{X,Y}$ be the subspace of the
path algebra $k\Quiver$ spanned by elements of the form 
\begin{align*}
\sum_{\{Z\in\IntersectionLattice: U_{i+1}\lessdot Z\lessdot U_{i-1}\}}
\Big(U_0 \arrow U_1 \arrows U_{i-1} \arrow Z \arrow U_{i+1}\arrows
U_{l-1}\arrow U_l\Big),
\end{align*}
where $0<i<l$, $U_0=X$ and $U_l=Y$. Note that $M_{X,Y}$ is a
subspace of $Y\ker(\varphi)X$ since
\begin{align*}
\sum_{\{Z: U_{i+1}\lessdot Z\lessdot U_{i-1}\}} (U_{i-1} \arrow Z
\arrow U_{i+1}) \in \ker(\varphi).
\end{align*}
Thus, $\dim(M_{X,Y})\leq\dim(Y\ker(\varphi)X)$. We show below that this is
an equality, which implies that $Y\mc I X = Y(\ker\varphi)X$, from
which it follows that $\ker\varphi=\mc I$.

We compute the dimension of the quotient space $Y(k\Quiver)X/M_{X,Y}$ using
results from \emph{poset cohomology} \cite{Wachs2007}. The poset
obtained by reversing the order on $\IntersectionLattice$ is a
geometric lattice and so the dual of the poset
$P=\{Z\in\IntersectionLattice : Y\leq Z \leq X\}$ is also a geometric
lattice \cite[Proposition 3.8]{Stanley2007}. The poset cohomology of
$P$ is isomorphic to the vector space $Y(k\Quiver)X/M_{X,Y}$ and its
dimension is known to be $|\mu(Y,X)|$, where $\mu$ is the M\"obius
function of $\IntersectionLattice$ \cite{Folkman1966,Bjorner1992}.
This is also the dimension of $e_Yk\Faces e_X$ because $\sum_{Y\leq
X}\dim(e_Yk\Faces e_X)=\dim(k\Faces e_X)$ counts the number of faces
of support $X$ (this follows from Lemma \ref{l:IdempotentLemma}; for a
proof see \cite[\S12]{Saliola2007a}) or \cite[Proposition
6.4]{Saliola2008a} and so does $\sum_{Y\leq X}|\mu(Y,X)|$
\cite{Zaslavsky1975}. Therefore,
$\dim(Y(k\Quiver)X/M_{X,Y})\leq\dim(e_Yk\Faces e_X)$. In particular,
$\dim(M_{X,Y})\geq\dim(Y(k\Quiver)X)-\dim(e_Yk\Faces e_X) =
\dim(Y(\ker\varphi)X)$. 
\end{proof}

\begin{Lemma}
\label{l:PhiIsWEquivariant}
$\varphi$ is $W$-equivariant.
\end{Lemma}

\begin{proof}
We need only show that $w(\varphi(P))=\varphi(w(P))$ for every path
$P$ in $\Quiver$ and every $w\in W$.

If $P$ is a path of length $0$, then $P$ is a vertex. Thus,
$w(\varphi(P)) = w(e_P) = e_{w(P)} = \varphi(w(P))$ for all $w \in W$.

If $P=(X\arrow Y)$ is an arrow in $\Quiver$, then for all $w\in W$,
\begin{align*}
  w(\varphi(X\arrow Y)) 
 &= w\Big(\ell_Y \Big([y:x] x + [y:x']x'\Big)e_X\Big) \\
 &= \ell_{w(Y)} \Big([y:x] w(x) + [y:x'] w(x')\Big) e_{w(X)}.
\end{align*}
It follows directly from Equation \eqref{e:IncidenceNumbers} and 
Equation \eqref{e:OrientationNumbers} that
$[x: y]$ is equal to 
 $\sigma_{\supp(x)}(w) \sigma_{\supp(y)}(w) [w(x):w(y)]$,
so
\begin{align*}
\Big([y:x] w(x) &+ [y:x'] w(x')\Big) \\
&= \sigma_Y(w)\sigma_X(w) \Big([w(y):w(x)] w(x) + [w(y):w(x')]
w(x')\Big).
\end{align*}
Hence, $w(\varphi(X\arrow Y)) = \varphi(w(X \arrow Y))$.

Since $w \Path Xp = w(X_{p-1}\arrow X_p)\cdots w (X_0\arrow X_1)$,
the result follows.
\end{proof}

This establishes Theorem \ref{t:WEquivariantPhi}. As an immediate
corollary, we get that $\Quiver$ is the quiver of $k\Faces$.

\begin{Corollary}
\label{c:QuiverOfSemigroupAlgebra}
$\Quiver$ is the quiver of $k\Faces$.
\end{Corollary}

\begin{proof}
From Theorem \ref{t:WEquivariantPhi}, $\varphi:k\Quiver\to k\Faces$ is
a surjective $k$-algebra homomorphism that satisfies
$0=F^{n+1}\subseteq\ker(\varphi)\subseteq F^2$, where $F$ is the ideal
in $k\Quiver$ generated by the arrows and $n=\dim(V)$. Therefore, by
Theorem \ref{t:ARS}, $\Quiver$ is the quiver of $k\Faces$.
\end{proof}

\section{On the Quiver of $\InvariantSubalgebra$}
\label{s:OnTheQuiverOfTheInvarviantSubalgebra}

Let $\InvariantQuiver$ denote the quiver of $\InvariantSubalgebra$.
This section explores some implications of Theorems \ref{t:CompSystem}
and \ref{t:WEquivariantPhi} for the structure of $\InvariantQuiver$.
Since $\DescentAlgebra$ is anti-isomorphic to $\InvariantSubalgebra$,
the quiver of $\DescentAlgebra$ is $\InvariantQuiver\opp$, the quiver
obtained from $\InvariantQuiver$ by reversing its arrows. So the
results below also apply to $\DescentAlgebra$ and
$\InvariantQuiver\opp$. In the next two sections we use Theorem
\ref{t:WEquivariantPhi} to compute the quiver of
$\InvariantSubalgebra[\Sn]$ and the quiver of
$\InvariantSubalgebra[\Bn]$.

Our first result deals with the vertices of $\InvariantQuiver$. Since
they correspond to idempotents in a \CompleteSystem for
$\InvariantSubalgebra$, Theorem \ref{t:CompSystem} implies that
$\InvariantQuiver$ has one vertex for each orbit $\Orbit \in
\IntersectionLattice/W$.

\begin{Proposition}
\label{p:VerticesOfInvariantQuiver}
$\InvariantQuiver$
has exactly one vertex for each $W$-orbit of elements in
$\IntersectionLattice$, where $\IntersectionLattice$ is the
intersection lattice of the reflection arrangement of $W$.
\end{Proposition}

Combined with Remark \ref{r:IntLatAndSubsetsOfS}, this implies that
the quiver of $\DescentAlgebra$ has exactly one vertex for each
equivalence class of subsets of $S$. 

The next observation will be the main tool in the remainder of this
section. It gives a sufficient condition for there to be no arrow
between $2$ given vertices in $\InvariantQuiver$.

\begin{Lemma}
\label{l:ConditionForNoArrow}
Let $\Orbit,\Orbit'\in\IntersectionLattice/W$ be vertices of
$\InvariantQuiver$. If for every path $P$ in $\Quiver$ that begins at
a vertex in $\Orbit'$ and ends at a vertex in $\Orbit$ there exists $w
\in W$ such that $w(P)=-P$, then there is no arrow from $\Orbit'$ to
$\Orbit$ in $\InvariantQuiver$.
\end{Lemma}

\begin{proof}
It follows from the definition of the quiver of an algebra
(\S\ref{s:QuiverOfSplitBasicAlgebra}) that if the vector space
$\varepsilon_\Orbit\InvariantSubalgebra\varepsilon_{\Orbit'}$ is the
zero vector space, then there is no arrow $\Orbit'\arrow\Orbit$. We'll
show that this vector space is zero if the hypothesis holds.

It follows from Theorem
\ref{t:WEquivariantPhi} that $\varphi$ restricts to a surjection 
$\nu_\Orbit \InvariantQuiverAlgebra \nu_{\Orbit'}\twoheadrightarrow
\varepsilon_\Orbit \InvariantSubalgebra \varepsilon_{\Orbit'}$,
where $\nu_\Orbit=\sum_{X\in\Orbit}X$ for each
$\Orbit\in\IntersectionLattice/W$.
We'll show $\nu_\Orbit\InvariantQuiverAlgebra\nu_{\Orbit'}=0$.
This subspace is spanned by elements of the form 
$\sum_{u\in W} u(P)$, where $P$ is a path of $\Quiver$ that begins at
a vertex in $\Orbit'$ and ends at a vertex in $\Orbit$. The hypothesis
states that $w(P)=-P$ for some $w\in W$, so
\begin{gather*}
\sum_{u\in W} u(P)
=\sum_{u\in W} u(w(P))
=-\left(\sum_{u\in W} u(P)\right).
\end{gather*}
Therefore, $\sum_{u\in W} u(P)=0$. So
$\nu_\Orbit\InvariantQuiverAlgebra\nu_{\Orbit'}=0$.
\end{proof}

Our first consequence of this lemma is that $\InvariantQuiver$
contains no oriented cycles. 
\begin{Proposition}
\label{p:NoOrientedCycles}
If $\Orbit' \arrow \Orbit$ is an arrow in $\InvariantQuiver$, then
$\Orbit < \Orbit'$ in $\IntersectionLattice/W$. In particular,
$\InvariantQuiver$ does not contain any oriented cycles.
\end{Proposition}
\begin{proof}
If $\Path Xl$ is a path in $\Quiver$, then $X_l\leq X_0$. In
particular, $\Orbit_{X_l} \leq \Orbit_{X_0}$. So if
$\Orbit\not<\Orbit'$, then the condition of Lemma
\ref{l:ConditionForNoArrow} is vacuously satisfied since there are no
paths in $\Quiver$ from a vertex in $\Orbit'$ to a vertex in $\Orbit$.
Therefore, there is no arrow from $\Orbit'$ to $\Orbit$ in
$\InvariantQuiver$. It follows that $\InvariantQuiver$ cannot contain
an oriented cycle.
\end{proof}

\begin{Corollary}
\label{c:DescentAlgebraIsQuasiHereditary}
The algebra $\InvariantSubalgebra$ is a quasi-hereditary algebra.
\end{Corollary}

This result follows from the definition of a quasi-hereditary algebra
since $\InvariantQuiver$ contains no oriented cycles (for an
introduction to quasi-hereditary algebras, see Vlastimil Dlab's
appendix to \cite{DrozdKirichenko1994}). Associated to every
quasi-hereditary algebra $A$ is a distinguished module $T$, called the
\emph{characteristic tilting module} of $A$, and the \emph{Ringel
dual} of $A$ is the algebra $\End_A(T)$; it develops that the Ringel
dual of $A$ is Morita equivalent to $A$ \cite{Ringel1991}. It would be
interesting to identify the characteristic tilting module and the
Ringel dual of $\InvariantSubalgebra$ and $\DescentAlgebra$.

Our next result shows that $\InvariantQuiver$ contains at least one
isolated vertex.

\begin{Proposition}
\label{p:NoArrowsFromTop}
There are no arrows in $\InvariantQuiver$ beginning at the vertex
$\Orbit_V$, where $V$ is the $W$-orbit of the ambient
vector space $V$ of the reflection arrangement of $W$.
\end{Proposition}
\begin{proof}
Let $\Path Xl$ be a path in $\Quiver$ with $X_0=V$. Let $w\in W$ denote
the reflection in the hyperplane $X_1$. Then 
\begin{align*}
w\Path Xl 
&= \sigma_{X_0}(w)\sigma_{X_l}(w)
\left( w(X_0)\arrow\adots\arrow w(X_l) \right) \\
&= -\Path{X}l
\end{align*}
since $w$ fixes pointwise all the subspaces $X_1,X_2,\ldots,X_l$
and changes the orientation of $X_0$.
By Lemma \ref{l:ConditionForNoArrow}, there is no arrow
in $\InvariantQuiver$ beginning at $\Orbit_V=\{V\}$.
\end{proof}

The poset $\IntersectionLattice/W$ is a ranked poset, with the rank of
an element $\Orbit\in\IntersectionLattice/W$ equal to the rank of any
$X\in\Orbit$ as an element of $\IntersectionLattice$ (which is
$\dim(X)-\dim(\cap_{H\in\Arrangement} H)$). As we will see in Theorem
\ref{t:QuiverOfDesAlgA}, if $W$ is the symmetric group $\Sn$,
$n\geq2$, then the existence of an arrow $\Orbit'\arrow\Orbit$ in the
quiver of $\DescentAlgebra$ implies that $\Orbit\lessdot\Orbit'$ in
$\IntersectionLattice/W$. The next result shows that this is not
necessarily true for other $W$.

\begin{Proposition}
\label{p:CoverRelationsNotGivingArrows}
If $W$ is a finite Coxeter group of type $A_1$, $B_n$, $D_{2n}$,
$E_7$, $E_8$, $F_4$, $I_2(2n)$, $H_3$ or $H_4$, then there is no arrow
in $\InvariantQuiver$ from $\Orbit'$ to $\Orbit$ if the difference
between their ranks in $\IntersectionLattice/W$ is odd. In particular,
if $\Orbit\lessdot\Orbit'$, then there is no arrow from $\Orbit'$ to
$\Orbit$ in $\InvariantQuiver$.
\end{Proposition}

\begin{proof}
If the type of $W$ is one of those listed above, then then $W$
contains the transformation $w(\vec v)=-\vec v$ for $\vec v \in V$
\cite[Lemma 27.2]{Kane2001}. Since $\sigma_X(w)=(-1)^{\dim(X)}$, we
have $\sigma_Y(w)\sigma_X(w)=-1$ if and only if $\dim(X)+\dim(Y)$ is
even.

If $\Path Xl$ is a path from $X_0 \in \Orbit'$ to $X_l \in \Orbit$,
then the hypothesis on the difference between the ranks of $\Orbit'$
and $\Orbit$ in $\IntersectionLattice/W$ implies that
$\dim(X_0)+\dim(X_l)$ is odd. Therefore,
\begin{align*}
w\Path Xl = \sigma_{X_0}(w)\sigma_{X_l}(w) \Path Xl =-\Path Xl.
\end{align*}
The result now follows from Lemma \ref{l:ConditionForNoArrow}.
\end{proof}

Combined with Proposition \ref{p:NoArrowsFromTop} the above result
implies the following.

\begin{Corollary}
If the type of $W$ is one of those listed in Proposition
\ref{p:CoverRelationsNotGivingArrows}, then $\InvariantQuiver$
contains at least three connected components.
\end{Corollary}

Recall that the \defn{Loewy length} of an algebra $A$ is the smallest
integer $\ell$ such that $\rad^\ell(A)=0$. If $W$ belongs to one of
the types listed in Proposition \ref{p:CoverRelationsNotGivingArrows},
then that result can be used to give an upper bound on the Loewy
length of $\InvariantSubalgebra$.

\begin{Proposition}
Let $(W,S)$ be a Coxeter system and let $n=|S|$. If $W$ is of type
$A_1$, $B_m$, $D_{2m}$, $E_7$, $E_8$, $F_4$, $I_2(2m)$, $H_3$ or
$H_4$, then the Loewy length of $\DescentAlgebra$ is at most
$\frac{n+1}2$.
\end{Proposition}
\begin{proof}
Since $\InvariantQuiver$ contains no oriented cycles, \emph{one plus}
the length of the longest path in $\InvariantQuiver$ is an upper bound
on the Loewy length of $\InvariantSubalgebra$. So we bound the length
of the longest path in $\InvariantQuiver$. Suppose
$\Orbit_0\arrow\Orbit_1\arrows\Orbit_l$ is a path in
$\InvariantQuiver$ with $l\geq1$. Since $n=|S|$ is the rank of the
poset $\IntersectionLattice$, Proposition \ref{p:NoArrowsFromTop}
implies that $\rank(\Orbit_0)\leq n-1$. Combined with 
Proposition \ref{p:CoverRelationsNotGivingArrows}, we obtain that
\begin{align*}
n-1\geq\Big(\rank(\Orbit_0)-\rank(\Orbit_l)\Big)
\geq\sum_{i=1}^l\Big(\rank(\Orbit_{i-1})-\rank(\Orbit_i)\Big)\geq 2l.
\end{align*}
Thus, $l\leq\frac{n-1}2$, and so the Loewy length of $\DescentAlgebra$
is at most $\frac{n+1}2$.
\end{proof}

These upper bounds are in fact equalities \cite{BonnafePfeiffer2008}.
This approach of bounding the length of the longest path in the quiver
was also used in \cite{Saliola2008b} to determine the Loewy length of
the descent algebra of type $D_{2m+1}$, the only case not covered by
earlier results \cite{BonnafePfeiffer2008}.

We have seen that the surjection
$\varphi:(k\Quiver)^W\twoheadrightarrow\InvariantSubalgebra$ plays an
important role in deducing information about $\InvariantQuiver$. The
next result, in conjunction with Theorem \ref{t:ARS}, explains why
this is the case. 

\begin{Theorem}
Suppose $\psi:k\InvariantQuiver\twoheadrightarrow\InvariantSubalgebra$
is a surjection of $k$-algebras with an admissible kernel and let
$\varphi:k\Quiver\twoheadrightarrow k\Faces$ denote the $k$-algebra
surjection of Theorem \ref{t:WEquivariantPhi}. Then $\psi$ factors
through $\varphi$.
\end{Theorem}

\begin{proof}
For each $\Orbit\in\IntersectionLattice/W$, let
$\nu_\Orbit=\sum_{X\in\Orbit}X\in k\Quiver$. The elements
$\varepsilon_\Orbit=\varphi(\nu_\Orbit)$ form a \CompleteSystem for
$\InvariantSubalgebra$ (Theorem \ref{t:CompSystem}), and so do the
elements $f_\Orbit=\psi(\Orbit)$. 
Since $\varepsilon_\Orbit$ and $f_\Orbit$ lift the same idempotent in
$\InvariantSubalgebra/\rad\InvariantSubalgebra
\cong(k\IntersectionLattice)^W$, there exists
$u_\Orbit\in\InvariantSubalgebra$ such that
$\varepsilon_\Orbit=u_\Orbit f_\Orbit u_\Orbit\inv$
\cite[Theorem 1.7.3]{Benson1998:I}. Let $u=\sum_{\Orbit}
\varepsilon_\Orbit u_\Orbit f_\Orbit$. Then $u\inv=\sum_{\Orbit}
f_\Orbit u_\Orbit\inv \varepsilon_\Orbit$ and
$\psi(\Orbit)=f_\Orbit=u\inv\varepsilon_\Orbit
u=u\inv\varphi(\nu_\Orbit)u$ for all
$\Orbit\in\IntersectionLattice/W$.

If $\Orbit'\arrow\Orbit$ is an arrow in $\InvariantQuiver$, then
$\psi(\Orbit'\arrow\Orbit)$ is a nonzero element of the subspace
$f_\Orbit\InvariantSubalgebra f_{\Orbit'}=
u\inv(\varepsilon_\Orbit\InvariantSubalgebra\varepsilon_{\Orbit'})u$.
Since $\varphi$ is surjective, there exists an element
$\rho_{(\Orbit'\arrow\Orbit)}$ in $\nu_\Orbit k\Quiver\nu_{\Orbit'}$
such that
$\psi(\Orbit'\arrow\Orbit)=u\inv\varphi(\rho_{(\Orbit'\arrow\Orbit)})u$,
and there exists $U\in k\Quiver$ such that $\varphi(U)=u$. Since
$\Quiver$ contains no oriented cycles, $\varphi(U)$ is invertible if
and only if $U$ is invertible. Thus, $U$ is invertible.
Let $\xi:k\InvariantQuiver\to k\Quiver$ be the homomorphism defined on
the vertices and arrows of $\InvariantQuiver$ by
$\xi(\Orbit)=U\inv\nu_\Orbit U$ and
$\xi(\Orbit'\arrow\Orbit)=U\inv\rho_{(\Orbit'\arrow\Orbit)}U$. It follows
that $\psi(P) = (\varphi\circ\xi)(P)$ for all $P\in
k\InvariantQuiver$. 
\end{proof}

\section{The Quiver of $\InvariantSubalgebra[\Sn]$}
\label{s:QuiverDesAlgA}

In this section we determine the quiver of
$\InvariantSubalgebra[\Sn]$. We begin by fixing notation. Throughout,
let $\Arrangement$ be the reflection arrangement of the symmetric
group $\Sn$, let $k\Faces$ and $\IntersectionLattice$ be the face
semigroup algebra and the intersection lattice of $\Arrangement$,
respectively, and let $\varphi:k\Quiver\to k\Faces$ be the map defined
in Theorem \ref{t:WEquivariantPhi}. Recall from
\S\ref{ss:SymmetricGroupExample} that an \defn{integer partition} of
$n\in \mathbb N$ is a collection of positive integers that sum to $n$.

\begin{Theorem}
\label{t:QuiverOfDesAlgA}
The quiver of $\InvariantSubalgebra[\Sn]$ is the directed graph with
one vertex $\nu_p$ for each integer partition $p$ of $n$ and exactly
one arrow $\nu_p \arrow \nu_q$ if and only if $q$ is obtained from $p$
by adding two distinct elements of $p$.
\end{Theorem}

\begin{figure}[!ht]
\ifnoxypic\xypicmessage\begin{comment}\fi
\begin{gather*}
\xymatrix@u@C=2em@R=2ex{
& 1111111\\
& 211111 \ar[dl] \\
31111 \ar[d]  & 22111 \ar[d] \\
4111 \ar[d]  & 3211 \ar[d] \ar[dl] \ar[dr] & 2221 \ar[dr] & \\
511 \ar[d]  & 421 \ar[dl] \ar[d] \ar[dr] & 331 \ar[d]
& 322 \ar[dll]  \\
61 \ar[dr] & 52 \ar[d] & 43 \ar[dl] \\
& 7}
\end{gather*}
\ifnoxypic\end{comment}
 \fi
\caption{The quiver of $\InvariantSubalgebra[{\Sn[7]}]$.}
\label{f:S7Quiver}
\end{figure}

\begin{proof}
Let $\InvariantQuiver$ be the quiver defined in the statement of the
theorem. We define a homomorphism of $k$-algebras $\psi:
k\InvariantQuiver\to\InvariantSubalgebra[\Sn]$ with an admissible
kernel. It then follows from Theorem \ref{t:ARS} that
$\InvariantQuiver$ is the quiver of $\InvariantSubalgebra[\Sn]$. See
\S\ref{ss:SymmetricGroupExample} for definitions.

\textbf{Definition of $\psi$.} For $X\in\IntersectionLattice$, write
$\pi(X)=\{B_1,\ldots,B_r\}$, where $|B_1|\geq\cdots\geq|B_r|$, for the
set partition associated to $X$, and let
$\rho(X)=(|B_1|,|B_2|,\ldots,|B_r|)$. Note that two elements $X$ and
$X'$ in $\IntersectionLattice$ are in the same $S_n$-orbit if and only
if $\rho(X)=\rho(X')$.

Define $\psi$ on the vertices $\nu_p$ of $\InvariantQuiver$ by
\begin{gather*}
\psi(\nu_p) = \sum_{X\in\IntersectionLattice\atop\rho(X)=p}
\varphi(X). 
\end{gather*}
If $\nu_p\arrow\nu_q$ is an arrow in $\InvariantQuiver$, then fix an
arrow $X\arrow Y$ in $\Quiver$ with $\rho(X)=p$ and $\rho(Y)=q$, and
define
\begin{align*}
\psi(\nu_p\arrow\nu_q) = \sum_{ w \in \Sn } \varphi(w(X \arrow Y)).
\end{align*}
We argue that $\psi$ extends to a unique $k$-algebra homomorphism
$\psi:\InvariantQuiver\to\InvariantSubalgebra$. By Theorem
\ref{t:ExtendingMapsToMorphisms}, we need to show that 
$\psi(\nu_p\arrow\nu_q)=\psi(\nu_q)\psi(\nu_p\arrow\nu_q)\psi(\nu_p)$
for all arrows $\nu_p\arrow\nu_q$ in $\InvariantQuiver$. Well,
\begin{align*}
\psi(\nu_p \arrow \nu_q) \psi(\nu_p)
&= \sum_{w \in \Sn} \varphi\left(w(X \arrow Y) \sum_{\rho(Z)=p} Z\right) \\
&= \sum_{w \in \Sn} \varphi\left(w(X \arrow Y) \sum_{\rho(Z)=p} w(Z)\right) \\
&= \sum_{w \in \Sn}(\varphi\circ w)\left((X \arrow Y) \sum_{\rho(Z)=p} Z\right) \\
&= \sum_{w \in \Sn}(\varphi\circ w)\left((X \arrow Y) X\right) \\
&= \sum_{w \in \Sn} \varphi\left(w\left( X \arrow Y \right) \right)
= \psi(\nu_p \arrow \nu_q).
\end{align*}
Similarly, $\psi(\nu_q)\psi(\nu_p\arrow\nu_q)=\psi(\nu_p\arrow\nu_q)$. 

\textbf{The kernel of $\psi$ is admissible.}
We next argue that the kernel of $\psi$ is an admissible ideal of
$k\InvariantQuiver$. Recall that an ideal of a path algebra is
admissible if every element in the ideal is a linear combination of
paths of length at least two. Suppose $a\in\ker(\psi)$. By multiplying
$a$ on the left and right by vertices of $\InvariantQuiver$, we can
suppose that $a$ is a linear combination of paths that begin at
$\nu_p$ and end at $\nu_q$. If $\nu_p=\nu_q$, then $a$ is a scalar
multiple of a vertex. This can't happen as $\psi(\nu_q)$ is nonzero
because it is part of a \CompleteSystem (Theorem \ref{t:CompSystem}).
If $\nu_p\to\nu_q$ is an arrow, then $\psi(\nu_p\to\nu_q)=\sum_w
w(X\arrow Y)$. This is zero if and only if there exists $w\in\Sn$ such
that $w(X\arrow Y)=-(X\arrow Y)$. We show this happens if and only if
$q=\rho(Y)$ is obtained from $p=\rho(X)$ by adding two equal parts of
$p$. Then we are done, since if such a $w$ exists, then
$\nu_p\to\nu_q$ is not an arrow of $\InvariantQuiver$. 

Let $\pi(X)=\{B_1,\ldots,B_r\}$ and suppose $|B_i|=p_i$ for all $1\leq
i\leq r$. Since $Y\lessdot X$, the set partition $\pi(Y)$ is obtained
from $\pi(X)$ by merging two blocks $B_i$ and $B_j$. By re-indexing we
can suppose $i=1$ and $j=2$. If $p_1=p_2$, then any permutation
$\omega\in\Sn$ that maps $B_1$ to $B_2$ and $B_2$ to $B_1$ while
fixing the other blocks of $\pi(X)$ will satisfy $\omega(X\arrow
Y)=-(X\arrow Y)$. Suppose instead that $p_1\neq p_2$. If
$\omega\in\Sn$ with $\omega(X)=X$ and $\omega(Y)=Y$, then $\omega$
permutes the blocks of $\pi(X)$ and the blocks of $\pi(Y)$. It follows
that $\omega(B_1)=B_2$ and $\omega(B_2)=B_2$ since $p_1\neq p_2$. Let
$x$ and $y$ be the set compositions $(B_1,B_2,B_3,\ldots,B_m)$ and
$(B_1\cup B_2,B_3,\ldots,B_m)$, respectively. Then, $y\omega(x)=x$. So
$\omega(x)$ and $x$ correspond to faces of support $X$ that lie on the
same side of $Y$. Since $\omega$ does not swap the two half spaces of
$X$ determined by $Y$, the discussion following Theorem
\ref{t:WEquivariantPhi} implies $\omega(X\arrow Y)=(X\arrow Y)$. 

Thus, $a$ is a linear combination of paths of length at least two,
so $\ker(\psi)$ is an admissible ideal of $k\InvariantQuiver$. 

\textbf{$\psi$ is surjective.}
We show that $\psi(k\InvariantQuiver)+\rad^2\InvariantSubalgebra[\Sn]
= \InvariantSubalgebra[\Sn]$; the result then follows from standard
ring theory: \emph{if $A$ is a $k$-algebra and $A'$ is a
$k$-subalgebra of $A$ such that $A' + \rad^2(A) = A$, then $A' = A$}
\cite[Proposition 1.2.8]{Benson1998:I}. To do this we will use the
following result of Manfred Schocker \cite[Theorem
$9.10$]{Schocker2006}: $\rad^2\InvariantSubalgebra[\Sn] = \rad^2
\left( k\Faces \right) \cap \InvariantSubalgebra[\Sn].$ (This can be
proved using results of this paper; such a proof is outlined in
Theorem \ref{t:SchockerLemma}.) 

Since $\varphi: k\Quiver \to k\Faces$ is surjective (Theorem
\ref{t:WEquivariantPhi}), it follows that the elements $\sum_{w\in
\Sn} w(\varphi(P))$, where $P$ is a path in $\Quiver$, span
$\InvariantSubalgebra[\Sn]$. Furthermore, $\rad^2(k\Faces)$ is spanned
by elements the $\varphi(P)$, where $P$ is of length at least two.
Thus, if $P$ has length at least two, then $\sum_{w} w(\varphi(P))$ is
in $\rad^2(k\Faces)\cap\InvariantSubalgebra[\Sn]
=\rad^2\InvariantSubalgebra[\Sn]$. If $P$ has length zero, then $P=X$
is a vertex and
\begin{gather*}
\sum_{w \in \Sn} w(\varphi(X)) 
= \varphi\left(\sum_{w\in\Sn}w(X)\right)
= \lambda
\varphi\left(\sum_{Y\in\IntersectionLattice\atop\rho(Y)=\rho(X)}
Y\right)
= \lambda \psi(\nu_{\rho(X)}),
\end{gather*}
where $\lambda=|\{w\in W: w(X)=X\}|$.

It remains to show that $\sum_w \varphi(w(P))\in\img(\psi)$ if $P$ is
an arrow. We first show that if $X\arrow Y$ and $X'\arrow Y'$ are two
arrows with $X$ and $X'$ in the same $\Sn$-orbit and $Y$ and $Y'$ in
the same $\Sn$-orbit, then there exists a permutation $u$ such that
$u(X'\arrow Y')=\pm(X\arrow Y)$. Let $\pi(X)=\{B_1,B_2,\ldots,B_r\}$
and $\pi(X')=\{B'_1,B'_2,\ldots,B'_r\}$. Since $X$ and $X'$ are in the
same orbit, there exists a permutation $w$ mapping $X'$ to $X$. So we
can assume that $X'=X$. Up to a re-indexing of the blocks, $B_1\cup
B_2$ is a block of $\pi(Y)$ and $B_3\cup B_4$ is a block of
$\pi(Y')$. Since $Y$ and $Y'$ are in the same orbit, it follows that
$|B_1|=|B_3|$ and $|B_2|=|B_4|$. Therefore, any permutation that swaps
$B_1$ with $B_3$ and $B_2$ with $B_4$ will map $X$ to $X$ and $Y'$ to
$Y$.

Let $X\arrow Y$ be an arrow in $\Quiver$. If there is a $w\in\Sn$ such
that $w(X\arrow Y)=-(X\arrow Y)$, then $\sum_w w({X\arrow Y})=0$.
Otherwise, $\sum_w \varphi(w({X\arrow
Y}))=\pm\psi(\Orbit_X\arrow\Orbit_Y)$ by the above.
\end{proof}

Since the descent algebra $\DescentAlgebra[\Sn]$ is isomorphic to the
opposite algebra of $\InvariantSubalgebra[\Sn]$ (Theorem
\ref{t:Bidigare}), its quiver is obtained by reversing the arrows in
Theorem \ref{t:QuiverOfDesAlgA}. This quiver, as a directed graph,
appears in the work of Adriano Garsia and Christophe Reutenauer
\cite{GarsiaReutenauer1989}; see especially \S5 and the figures
contained therein. Manfred Schocker \cite[Theorem 5.1]{Schocker2004}
was the first to show that this is the quiver of
$\DescentAlgebra[\Sn]$ by using results of Dieter Blessenohl and
Hartmut Laue \cite{BlessenohlLaue1996,BlessenohlLaue2002}.

We also remark that the argument presented above can be used to find
the quiver of $\InvariantSubalgebra$ for arbitrary finite Coxeter
groups $W$ once the relationship between
$\rad^2(\InvariantSubalgebra)$ and
$\rad^p(k\Faces)\cap\InvariantSubalgebra$ is understood. We do this in
\S\ref{s:QuiverDesAlgB} for the finite Coxeter group of type $B$.

\subsection{Descending Loewy series of $\InvariantSubalgebra[\Sn]$}
\label{ss:LoewySeriesTypeA}

The proof of Theorem \ref{t:QuiverOfDesAlgA} relied on the case $m=2$
of the following result of Manfred Schocker.

\begin{Theorem}[Theorem $9.10$ of \cite{Schocker2006}] 
\label{t:SchockerLemma}
Let $k\Faces$ be the face semigroup algebra of the reflection
arrangement of the symmetric group $\Sn$. For all $m\in\mathbb N$, 
\begin{gather*}
\rad^m (k\Faces)^{\Sn} = \rad^m( k\Faces ) \cap
(k\Faces)^{\Sn}.
\end{gather*}
\end{Theorem}

Manfred Schocker proved this by constructing a basis of
$\rad^m(k\Faces)\cap\InvariantSubalgebra[\Sn]$ and noting that the
basis coincides, under an anti-isomorphism
$\InvariantSubalgebra[\Sn]\cong\DescentAlgebra[\Sn]$, to a basis of
$\rad^m\DescentAlgebra[\Sn]$ constructed by Dieter Blessenohl and
Hartmut Laue \cite{BlessenohlLaue1996}. 

This result can also be proved using just the theory developed in this
paper. In fact, in Theorem \ref{t:TypeBLoewySeries}, we prove the
corresponding result for the hyperoctahedral group, which is new. That
proof can be adapted to prove the above. We provide a very rough
outline of the argument and leave the details to the interested
reader.

\begin{proof}[Outline of a proof of Theorem \ref{t:SchockerLemma}]
To prove this, a lemma corresponding to Lemma
\ref{l:ExpansionLemmaTypeB}---and proved by arguing in the same
way---is needed:
\begin{align*}
\alpha\Norm\Big(X_0 \arrow\ccdots\arrow X_m\Big) 
&= 
\Norm\Big(X_1 \arrows X_m\Big) \Norm\Big(X_0 \arrow X_1\Big) \\
&\qquad - 
 \smashoperator{\sum_{t \in {\Sn}, t(X_1) = X_1, \atop 
  t(A \cup B) \neq A \cup B }} 
\quad \sigma_{X_1}(t)\sigma_{X_m}(t) 
 \Norm\Big( X_0 \arrow t(X_1) \arrows t(X_m)\Big),
\end{align*}
where $\Path Xm$ is a path in $\Quiver$ of length $m\geq2$, the set
partition $\Partition1$ is obtained by merging two blocks $A$ and $B$
of $\Partition0$, and $\alpha$ is the number of permutations that fix
$X_1$ and $A\cup B$.

Begin by reducing to the case $m=2$ by mimicking the proof of Theorem
\ref{t:TypeBLoewySeries}. To prove the case $m=2$, first establish the
containment
$\RAD[2]\subseteq\rad^2(k\Faces)\cap\InvariantSubalgebra[\Sn]$. For
the reverse containment, argue by contradiction: suppose that there
exists a path $P=\Path Xm$ in $\Quiver$ of length at least two such
that $\iNorm P\not\in\RAD[2]$; and of all such paths (that begin at
$X_0$), pick $P$ such that $|A|+|B|$ is maximal, where $A$ and $B$ are
the blocks of $\Partition0$ that are merged to get $\Partition1$. Then
argue as in Step 3 of the proof of Theorem \ref{t:TypeBLoewySeries}
that $A\cup B$ is not a block of $\Partition m$. This means that
$A\cup B$ is merged with some other block at some point. Argue as in
Step 2 of the proof to show (using the relations in the partition
lattice) that we can suppose that $\Partition2$ is obtained by merging
two blocks $C$ and $D$, where $|C|=|A|+|B|$. Then derive a
contradiction as in Step 4 of the proof, by examining the three cases:
$C,D\neq A\cup B$; $D=A\cup B\neq C$; $C=A\cup B\neq D$.
\end{proof}

\section{The Quiver of $\InvariantSubalgebra[\Bn]$}
\label{s:QuiverDesAlgB}

In this section we determine the quiver of
$\InvariantSubalgebra[\Bn]$. Throughout, let $\Arrangement$ be the
reflection arrangement of the hyperoctahedral group $\Bn$ (defined
below), let $k\Faces$ and $\IntersectionLattice$ be the face semigroup
algebra and the intersection lattice of $\Arrangement$, respectively,
and let $\varphi:k\Quiver\to k\Faces$ be the map defined in Theorem
\ref{t:WEquivariantPhi}. 

\subsection{The Coxeter group of type $\bs B$}
\label{ss:TypeBCoxGrp}
Let $n\in\mb N$. The \defn{Coxeter group of type B} and \defn{rank}
$n$, denoted by $\Bn$, is the finite group of orthogonal
transformations of $\mb R^n$ generated by reflections in the
hyperplanes 
\begin{align*}
\{\vec v \in \mb R^n : v_i = 0\},\quad
\{\vec v \in \mb R^n : v_i = v_j\},\quad
\{\vec v \in \mb R^n : v_i = -v_j\},
\end{align*}
where $i,j\in\{1,2,\ldots,n\}$ and $i\neq j$. This set of hyperplanes
is the \defn{reflection arrangement} of $\Bn$. We identify $\Bn$ with
the group of \emph{signed permutations} as follows.
For $n\in\mb N$, let $[n] = \{1,2,\ldots,n\}$ and let $[\pm n]=[n]
\cup (-[n])$. A \defn{signed permutation} of $[\pm n]$ is a
permutation $w$ of the set $[\pm n]$ satisfying $w(-i)=-w(i)$ for all
$i\in[n]$. 
Every signed permutation $w$ induces an orthogonal
transformation of $\mb R^n$ by permuting and negating coordinates.
Moreover, any transformation in $\Bn$ arises in this fashion.

For any $A\subseteq[\pm n]$ let $\overline A=\{-i: i\in A\}$.
Under the above identification the intersection lattice of the type
$B$ arrangement is identified with the sublattice $\PartitionLatticeB$
of set partitions of $[\pm n]$ of the form $\{B_1, \ldots, B_r, Z,
\overline B_r, \ldots, \overline B_1\}$, and where $Z$ can be empty
and satisfies $\overline Z=Z$ \cite[Theorem~4.1]{BarceloIhrig1999}.

To simplify notation, we let $\pi(X)$ denote the set partition of
$[\pm n]$ induced by $X\in\IntersectionLattice$, and we let
$\{B_1,\ldots,B_r;Z\}$ denote the set partition $\{B_1, \ldots, B_r,
Z, \overline B_r, \ldots, \overline B_1\}$. The set $Z$ is called the
\defn{zero block} and the other sets are called \defn{nonzero blocks}.
Under this isomorphism the action of $\Bn$ on $X \in
\IntersectionLattice$ is given by permuting the elements of $\pi(X)$.
That is, $\pi(w(X)) = w(\pi(X))$ for all $w\in \Bn$ and
$X\in\IntersectionLattice$.

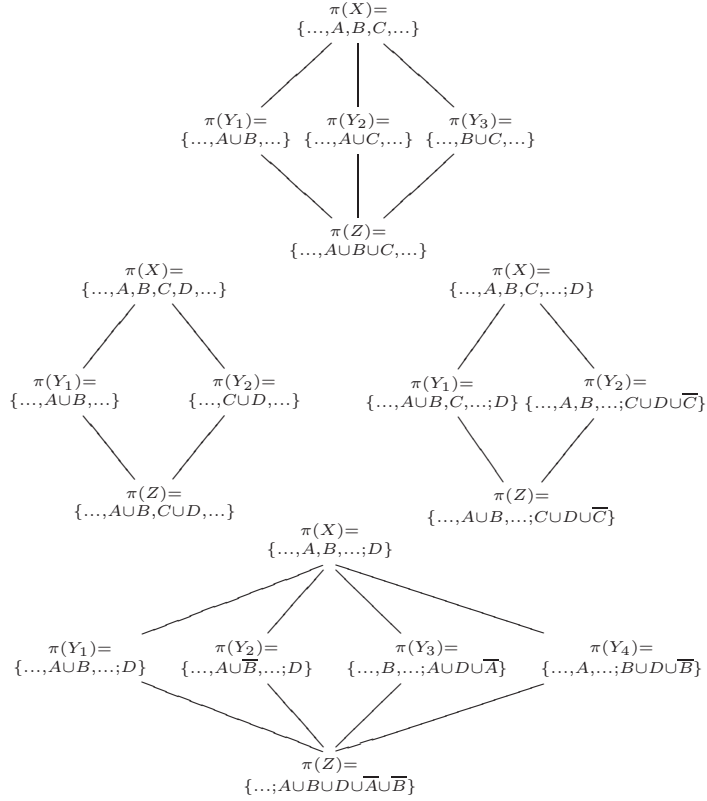
\begin{figure}
\ifhrule\hrule\fi
{\footnotesize 
\begin{gather*}
\xymatrix@C=-1em{
& \pi(X)\mathord=\atop\{\cldots,A,B,C,\cldots\} \ar@{-}[dl] \ar@{-}[d] \ar@{-}[dr] \\
\pi(Y_1)\mathord=\atop\{\cldots,A \cup B,\cldots\} \ar@{-}[dr] 
&\pi(Y_2)\mathord=\atop\{\cldots,A \cup C,\cldots\} \ar@{-}[d] 
&\pi(Y_3)\mathord=\atop\{\cldots,B \cup C,\cldots\} \ar@{-}[dl] \\
&\pi(Z)\mathord=\atop\{\cldots,A \cup B \cup C,\cldots\} 
}
\\[-2ex]
\xymatrix@C=-3em{
& \pi(X)\mathord=\atop\{\cldots,A,B,C,D,\cldots\} \ar@{-}[dl] \ar@{-}[dr] \\
\pi(Y_1)\mathord=\atop\{\cldots,A \cup B,\cldots\} \ar@{-}[dr] 
& 
&\pi(Y_2)\mathord=\atop\{\cldots,C \cup D,\cldots\}\ar@{-}[dl] \\
&\pi(Z)\mathord=\atop\{\cldots,A \cup B, C \cup D,\cldots\}}
\qquad
\xymatrix@C=-5em{
& \pi(X)\mathord=\atop\{\cldots,A,B,C,\cldots;D\} \ar@{-}[dl] \ar@{-}[dr] \\
\pi(Y_1)\mathord=\atop\{\cldots,A \cup B,C,\cldots;D\} \ar@{-}[dr] 
& 
&\pi(Y_2)\mathord=\atop\{\cldots,A,B,\cldots;C\cup D\cup \overline{C}\}\ar@{-}[dl] \\
&\pi(Z)\mathord=\atop\{\cldots,A\cup B,\cldots;C\cup D\cup \overline{C}\}} 
\\[-2ex]
\xymatrix@C=0em{
&&&\save[]+<0em,1em>*{%
 \pi(X)\mathord=\atop\{\cldots,A,B,\cldots;D\}}\restore\ar@{-}[drrr]\ar@{-}[dr]\ar@{-}[dl]\ar@{-}[dlll]\\
\pi(Y_1)\mathord=\atop{\{\cldots,A\cup B,\cldots;D\}}&&
\pi(Y_2)\mathord=\atop{\{\cldots,A\cup \overline{B},\cldots;D\}}&&
\pi(Y_3)\mathord=\atop{\{\cldots,B,\cldots;A\cup D\cup\overline{A}\}}&&
\pi(Y_4)\mathord=\atop{\{\cldots,A,\cldots;B\cup D\cup\overline{B}\}}\\
&&&\save[]+<0em,-1em>*{%
\pi(Z)\mathord=\atop\{\cldots;A\cup B\cup D\cup\overline{A}\cup\overline{B}\}}\restore
 \ar@{-}[urrr]\ar@{-}[ur]\ar@{-}[ul]\ar@{-}[ulll]
}
\end{gather*}}
\caption{The four types of intervals of length two in the lattice
$\PartitionLatticeB$ of set partitions of type B.}
\label{f:PartitionLatticeTypeBTwoIntervals}
\ifhrule\hrule\fi
\end{figure}

The intervals of length two in $\PartitionLatticeB$ play an important
role in what follows. So we quickly describe them. If $P' \lessdot P$
is a cover relation in $\PartitionLatticeB$, then either $P'$ is
obtained from $P$ by merging two distinct nonzero blocks of $P$, or
$P'$ is obtained from $P$ by merging a nonzero block $B$ with
$\overline B$ and the zero block of $P$. It follows that there are
four types of intervals of length two in $\PartitionLatticeB$, which
are illustrated in Figure \ref{f:PartitionLatticeTypeBTwoIntervals}.

\subsection{The Quiver of $\InvariantSubalgebra[\Bn]$}
We now describe the quiver of $\InvariantSubalgebra[\Bn]$.

\begin{Theorem}
\label{t:QuiverOfDesAlgB}
The quiver of $\InvariantSubalgebra[\Bn]$ contains one vertex $\nu_p$
for each integer partition $p$ of $0,1,\ldots,n$, and $m_{p,q}$ arrows
from $\nu_p$ to $\nu_q$, where
\begin{align*}
m_{p,q} = 
\begin{cases}
2, &\text{if $q$ is obtained by adding $3$ distinct parts of $p$,}\\
1, &\text{if $q$ is obtained by adding $3$ parts of $p$, $2$ of which
are distinct,} \\
1, &\text{if $q$ is obtained by deleting $2$ distinct parts of $p$,} \\
0, &\text{otherwise.} \\
\end{cases}
\end{align*}
\end{Theorem}

The quiver of $\InvariantSubalgebra[{\Bn[6]}]$ is illustrated in
Figure~\ref{f:B6Quiver}.
\begin{figure}
\begin{gather*}
\xymatrix@u@C=1.5em@R=1.5ex{
& & & 111111 \\
& & 
 21111 \ar[dll] \ar[drrrr] & & 11111 
\\
 411 \ar[dd] \ar[ddrrrrrr]
 & 321 \ar@/^1ex/@<0.25ex>[ddl] \ar@/_1ex/@<0.25ex>[ddl] 
       \ar@/^2ex/[ddrrr] \ar@/^1.5ex/[ddrrrr] \ar@/^1.5ex/[ddrrrrr]
 & 222 
 & 311 \ar[ddll] \ar@/^/[ddrrr]
 & 221 \ar[ddlll] \ar[ddr]
 & 211 \ar[ddlll] \ar[ddr]
 & 111 
\\ \\
6 & 5 & 4 & & 3 & 2 & 1
}
\qquad
\xymatrix@u@C=0.75em@R=6ex{
& 
 & 3111 \ar[dll] \ar[drrrrrr]
 & 
 & 2211 \ar[dlll] \ar[drrr] \ar[dllll]
 & 2111 \ar[dll] \ar[drrr]
 & 1111 
 \\
51 & 42 & 33 & 41 & 32 & 31 & 22 & 21 & 11 \\
& & & & 
\emptyset
 \ar@{<-}@<0.5ex>[ullll]
 \ar@{<-}[ulll]
 \ar@{<-}[ul]
 \ar@{<-}[u]
 \ar@{<-}[ur]
 \ar@{<-}[urrr]
}
\end{gather*}
\caption{The quiver of $\InvariantSubalgebra[{\Bn[6]}]$.}
\label{f:B6Quiver}
\end{figure}

\begin{proof}
Let $\InvariantQuiver$ be the quiver with one vertex $\nu_p$ for each
integer partition $p$ of $0,1\ldots,n$, and $m_{p,q}$ arrows from the
vertex $\nu_p$ to $\nu_q$. We will use Theorem \ref{t:ARS} to show
that $\InvariantQuiver$ is the quiver of $\InvariantSubalgebraB$ by
constructing a surjective $k$-algebra morphism $\psi:
k\InvariantQuiver \twoheadrightarrow \InvariantSubalgebraB$ that has
an admissible kernel.

For each $X\in\IntersectionLattice$, let $\pi(X)=\{A_1,\ldots,A_r;Z\}$
with $|A_1|\geq\cdots\geq|A_r|$, and let $\rho(X)$ be the integer
partition $(|A_1|,|A_2|,\ldots,|A_r|)$. It follows that $X$ and $Y$
are in the same $\Bn$-orbit if and only if $\rho(X)=\rho(Y)$. 

\medskip
\textbf{Definition of $\psi$ on vertices.}
Let $\varphi:k\Quiver\to k\Faces$ denote the $k$-algebra homomorphism
of Theorem \ref{t:WEquivariantPhi}. Define a function $\psi$ on the
vertices of $\InvariantQuiver$ by
\begin{gather*}
\psi(\nu_p) = \sum_{\rho(X)=p} \varphi(X),
\end{gather*}
where $p$ is an integer partition of some $m\in\{0,1,\ldots,n\}$.

\medskip
\textbf{Definition of $\psi$ on arrows.} We define $\psi$ on the three
types of arrows of $\InvariantQuiver$ individually. See Figure
\ref{f:PartitionLatticeTypeBTwoIntervals} for the different types of
intervals of length two in $\PartitionLatticeB$.

Suppose $q$ is obtained by adding three distinct parts $p_1$, $p_2$ and
$p_3$ of $p$, where $p_1 > p_2 > p_3$, and let $\alpha_1^{p,q}$ and
$\alpha_2^{p,q}$ be the two arrows in $\InvariantQuiver$ from $\nu_p$
to $\nu_q$. Let $X\in\IntersectionLattice$ with $\rho(X)=p$, and let
$A$, $B$ and $C$ be three distinct blocks of $\pi(X)$ with
$|A|=p_1$, $|B|=p_2$ and $|C|=p_3$. Let $\pi(Y_1)$ be the set
partition obtained from $\pi(X)$ by merging $A$ and $B$, let
$\pi(Y_2)$ be the set partition obtained from $\pi(X)$ by merging
$A$ with $C$ and let $\pi(Z)$ be the set partition obtained from
$\pi(X)$ by merging $A$, $B$ and $C$. For $i\in\{1,2\}$, define
\begin{align*}
\psi(\alpha_i^{p,q}) =
\sum_{w\in \Bn} \varphi(w(X\arrow Y_i\arrow Z)).
\end{align*}

Suppose $q$ is obtained by adding three parts $p_1,p_2$ and $p_3$ of
$p$ with $p_1\neq p_2=p_3$, and let $\beta_{p,q}$ be the arrow in
$\InvariantQuiver$ from $\nu_p$ to $\nu_q$. Let $X,Y_1,Y_2$ and $Z$ be
as above. Define
\begin{align*}
\psi(\beta_{p,q}) =
\sum_{w\in \Bn} \varphi(w(X\arrow Y_1\arrow Z)).
\end{align*}

Finally, suppose that $q$ is obtained by deleting two distinct
parts $p_1$ and $p_2$ of $p$, and let $\gamma_{p,q}$ be the arrow in
$\InvariantQuiver$ from $\nu_p$ to $\nu_q$. Let
$X\in\IntersectionLattice$ with $\rho(X)=p$, and let $A$ and $B$
be two distinct blocks of $\pi(X)$ with $|A|=p_1$ and $|B|=p_2$.
Let
$\pi(Y_1)$ be the set partition obtained from $\pi(X)$ by merging $A$
and $B$, and let $\pi(Z)$ be the set partition obtained from
$\pi(X)$ by merging $A$, $B$, $\overline{A}$, $\overline{B}$
and the zero block of $\pi(X)$. Define
\begin{align*}
\psi(\gamma_{p,q}) =
\sum_{w\in \Bn} \varphi(w(X\arrow Y_1\arrow Z)).
\end{align*}

\textbf{Extension of $\psi$ to an algebra homomorphism.}
By Theorem \ref{t:ExtendingMapsToMorphisms}, $\psi$ extends to a
unique $k$-algebra homomorphism
$\psi:k\InvariantQuiver\to\InvariantSubalgebraB$ if the elements
$\psi(\nu_p)$ form a \CompleteSystem and if
$\psi(\nu_q)\psi(\nu_p\arrow\nu_q)\psi(\nu_p)=\psi(\nu_p\arrow\nu_q)$
for every arrow $\nu_p\arrow\nu_q$ in $\InvariantQuiver$. The first
condition follows from Theorem \ref{t:WEquivariantPhi}. Write $\psi(\nu_p
\arrow \nu_q)=\sum_w w(X\arrow Y\arrow Z)$ and note that
\begin{align*}
\psi(\nu_p \arrow \nu_q) \psi(\nu_p)
&=\sum_{w \in \Bn} 
  \varphi\left(w(X\arrow Y\arrow Z) \sum_{\rho(X')=p} X'\right)
  \\
&= \sum_{w \in \Bn} \varphi\left(w(X\arrow Y\arrow Z)
  \sum_{\rho(X')=p} w(X') \right)\\
&= \sum_{w \in \Bn} \varphi\left(w\left( \sum_{\rho(X')=p}
(X\arrow Y\arrow Z) X' \right) \right)\\
&= \sum_{w \in \Bn}\varphi\Big(w\left(X\arrow Y\arrow Z\right)\Big)
= \psi(\nu_p \arrow \nu_q).
\end{align*}
Similarly, $\psi(\nu_q)\psi(\nu_p\arrow\nu_q)=\psi(\nu_p\arrow\nu_q)$.

\medskip
\textbf{$\psi$ is surjective.}
Next we prove that $\psi:k\InvariantQuiver\to\InvariantSubalgebraB$ is
surjective. We show that $\psi(k\InvariantQuiver) +
\rad^2(\InvariantSubalgebraB) = \InvariantSubalgebraB$; the result
then follows from standard ring theory: \emph{if $A$ is a $k$-algebra
and $A'$ is a $k$-subalgebra of $A$ such that $A' + \rad^2(A) = A$,
then $A' = A$} \cite[Proposition 1.2.8]{Benson1998:I}. 
In order to do this we'll use a fact whose proof we defer to later
(Theorem \ref{t:TypeBLoewySeries}): that
$\rad^2(\InvariantSubalgebraB)=
\rad^4\left(k\Faces\right)\cap\InvariantSubalgebraB.$

Since $\varphi:k\Quiver\to k\Faces$ is surjective (Theorem
\ref{t:WEquivariantPhi}), the images of the paths $P$ of $\Quiver$
span $k\Faces$. It follows that the elements $\iNorm P$, form a
spanning set for $\InvariantSubalgebraB$ (recall that
$\Norm(P)=\sum_{w\in \Bn} w(P)$). Furthermore, $\rad^4(k\Faces)$ is
spanned by elements $\iNorm P$, where $P$ is of length at least four.
So if $P$ has length at least $4$, then $\iNorm P$ is in
$\rad^4(k\Faces)\cap\InvariantSubalgebraB$, so it is in
$\rad^2(\InvariantSubalgebraB)$ (Theorem \ref{t:TypeBLoewySeries}). It
remains to prove that $\iNorm P\in\img(\psi)$ if the length of $P$ is
less than $4$. If the length of $P$ is odd, then the signed
permutation $i\mapsto -i$ for all $i\in[\pm n]$ maps $P$ to $-P$, so
$\Norm(P)=0$. It remains to prove this for vertices and for paths of
length $2$.

Suppose $P=X$ is a vertex of $\Quiver$. If
$\lambda=|\{w\in\Bn:w(X)=X\}|$, then
\begin{align*}
\iNorm X=
\sum_{w\in \Bn} w(\varphi(X))=
\lambda \sum_{\rho(Y)=\rho(X)}\varphi(Y)=
\lambda \psi(\nu_{\rho(X)}).
\end{align*}
So $\iNorm X\in\img(\psi)$.

Suppose that $P=(X\arrow Y\arrow Z)$ is path of length two in
$\Quiver$. Let $p=\rho(X)$ and $q=\rho(Z)$. There are four cases to
consider, corresponding to the four types of intervals illustrated in
Figure \ref{f:PartitionLatticeTypeBTwoIntervals}.

\smallskip
\emph{Case 1. Suppose $q$ is obtained from $p$ by adding $p_1$ to
$p_2$ and deleting $p_3$, where $p_1$, $p_2$ and $p_3$ are three parts
of $p$.}
Since $\rho(Z)=q$, a nonzero block $C$ of $\pi(X)$ is contained in
the zero block of $\pi(Z)$. The signed permutation that negates the
elements of $C$ maps $P$ to $-P$, and so $\iNorm{P}=0\in\img(\psi)$.

\smallskip
\emph{Case 2. Suppose $q$ is obtained from $p$ by deleting $p_1$
and $p_2$, where $p_1$ and $p_2$ are two parts of $p$.}
Then there are two nonzero blocks $A$ and $B$ of $\pi(X)$, of sizes $p_1$
and $p_2$, respectively, that are contained in the zero block of
$\pi(Z)$. We'll show that if $P'=(X'\arrow Y'\arrow Z')$ is a
path with $\rho(X')=p$ and $\rho(Z')=q$, then $P'$ is in the
$\Bn$-orbit of a path from $X$ to $Z$. If $\rho(X')=p$, then there
exists $w\in \Bn$ such that $w(P')$ begins at $X$. Since $w(\pi(Z'))$
is obtained from $\pi(X)$ by merging two blocks $A'$ and $B'$ of
sizes $p_1$ and $p_2$, respectively, with the zero block of $\pi(X)$,
it follows that the signed permutation $u$ that swaps $A'$ with
$A$ and $B'$ with $B$ maps $w(P')$ to a path that begins at $X$
and ends at $Z$.

There are exactly four paths $P_i=(X\arrow Y_i\arrow Z)$, where
$i\in\{1,2,3,4\}$, in $\Quiver$ that begin at $X$ and end at $Z$: 
$\pi(Y_1)$ contains the block $A\cup B$;
$\pi(Y_2)$ contains the block $A\cup\overline{B}$;
the zero block of $\pi(Y_3)$ contains $A$;
the zero block of $\pi(Y_4)$ contains $B$.
The signed permutation that negates $A$ maps $P_3$ to $-P_3$, so 
$\Norm(P_3)=0$. Similarly, $\Norm(P_4)=0$. Since
$P_1+P_2+P_3+P_4\in\ker(\varphi)$ (Lemma \ref{l:KernelOfPhi}), we
have $\iNorm{P_1}=-\iNorm{P_2}$.

If $\iNorm{P_1}=0$, then $\Norm(P_1)\in\ker(\varphi)$, so $\Norm(P_1)$
is a scalar multiple of $P_1+P_2+P_3+P_4$. Since $P_3$ is not in the
orbit of $P_1$, it follows that $\Norm(P_1)=0$. So there exists a
signed permutation that maps $P_1$ to $-P_1$. This happens if and only
if $|A|=|B|$. So if $p_1=p_2$, then $\iNorm{P_1}=0$, and
if $p_1\neq p_2$, then $\iNorm{P_1}=\pm \psi(\gamma_{p,q})$. Since
$P\in\{P_1,P_2,P_3,P_4\}$, it follows that $\iNorm{P}\in\img(\psi)$.

\smallskip
\emph{Case 3. Suppose $q$ is obtained by adding three parts
$p_1$, $p_2$ and $p_3$ of $p$.}
Since $p$ and $q$ are partitions of the same integer, the zero blocks
of $\pi(X)$ and $\pi(Z)$ are the same. So there are two possibilities
for $\pi(Z)$: either $\pi(Z)$ contains the nonzero blocks $A\cup
B$ and $C\cup D$, or $\pi(Z)$ contains the nonzero block
$A\cup B\cup C$, where $A,B,C$ and $D$ are (nonzero)
blocks of $\pi(X)$. In the first case the signed permutation that
negates the elements of $A\cup B$ maps $P$ to $-P$, thus
$\Norm(P)=0$. 

So suppose $\pi(Z)$ contains the nonzero block $A\cup B\cup C$,
and that $|A|=p_1$, $|B|=p_2$ and $|C|=p_3$. Let $P'=(X'\arrow
Y'\arrow Z')$ be another path in $\Quiver$ with $\rho(X')=p$ and
$\rho(Z')=q$. Then either $\pi(Z')$ contains the nonzero block $A'\cup
B'\cup C'$, where $A',B',C'$ are blocks of $\pi(X')$
with $|A'|=|A|$, $|B'|=|B|$ and $|C'|=|C|$, or $\Norm(P')=0$ (as
above). In the former situation we have, by arguing as in Case 2, that
$P'$ is in the $\Bn$-orbit of a path that begins at $X$ and ends at
$Z$. This implies that $\psi(\nu_p\arrow\nu_q)=\pm \iNorm{P'}$ for
some path $P'$ beginning at $X$ and ending at $Z$.

Let $P_i=(X\arrow Y_i\arrow Z)$ for $i\in\{1,2,3\}$ be the three paths
in $\Quiver$ from $X$ to $Z$, where $A\cup B$ is a block of
$\pi(Y_1)$, $A\cup C$ is a block of $\pi(Y_2)$, and $B\cup C$ is a
block of $\pi(Y_3)$. 
If $p_1 > p_2 > p_3$, then the previous paragraph implies that
$\psi(\alpha_i^{p,q})=\pm\iNorm{P_i}$ for $i\in\{1,2\}$. Hence,
$\iNorm{P_1}$ and $\iNorm{P_2}$ are in $\img(\psi)$, and so
$\iNorm{P_3}\in\img(\psi)$ since $P_1+P_2+P_3\in\ker(\varphi)$
(Theorem \ref{t:WEquivariantPhi}). Therefore, $\iNorm P\in\img(\psi)$
since $P\in\{P_1,P_2,P_3\}$.

Suppose $p_1=p_2\neq p_3$ and suppose $\pi(Y_1)$ contains the block
$A\cup B$. It follows that
$\psi(\beta_{p,q})=\pm\iNorm{P_i}$ for some $i\in\{2,3\}$. The
signed permutation that swaps $A$ and $B$ maps $P_1$ to $-P_1$ and
so $\Norm(P_1)=0$. Since $P_1+P_2+P_3\in\ker(\varphi)$ (Lemma
\ref{l:KernelOfPhi}), it follows
that $\iNorm{P_2}=-\iNorm{P_3}=\pm \psi(\beta_{p,q})$. Thus,
$\iNorm P\in\img(\psi)$.

If $p_1=p_2=p_3$, then the argument in the previous paragraph implies
that $\Norm(P_i)=0$ for all $i\in\{1,2,3\}$. Hence, $\iNorm
P\in\img(\psi)$.

\smallskip
\emph{Case 4. Suppose that $q$ is obtained from $p$ by adding $p_1$ to
$p_2$ and by adding $p_3$ to $p_4$, where $p_1$, $p_2$, $p_3$ and
$p_4$ are parts of $p$.} 
If $q$ can also be obtained from $p$ by merging three parts of $p$,
then we can apply the argument of the previous case.
On the other hand, suppose $q$ is not obtained from $p$ by merging
three parts of $p$. Then $\pi(Z)$ contains the nonzero blocks $A\cup
B$ and $C\cup D$, where $A$, $B$, $C$, and $D$ are
(nonzero) blocks of $\pi(X)$ and $|A|=p_1$, $|B|=p_2$, $|C|=p_3$
and $|D|=p_4$.
The signed permutation that negates the elements of $A\cup B$ maps
$P$ to $-P$, so $\Norm(P)=0$. Hence, $\iNorm P\in\img(\psi)$.

\medskip
\textbf{The kernel of $\psi$ is admissible.}
To complete the proof we need to show that $\ker(\psi)$ is an
admissible ideal of $k\InvariantQuiver$. Recall that an ideal of a
path algebra is admissible if every element in the ideal is a linear
combination of paths of length at least two. Suppose $a\in\ker(\psi)$.
By multiplying $a$ on the left and right by vertices of
$\InvariantQuiver$, we can suppose that $a$ is a linear combination of
paths that begin at $\nu_p$ and end at $\nu_q$. If $\nu_p=\nu_q$, then
$a$ is a scalar multiple of a vertex. This implies $a=0$ because
$\psi(\nu_q)$ is nonzero: it belongs to a \CompleteSystem (Theorem
\ref{t:WEquivariantPhi}). If $a$ is a linear combination of arrows
that begin at $\nu_p$ and end at $\nu_q$, then there are three cases
to consider depending on the type of the arrows. We will show that
$\psi(\alpha_1^{p,q})$ and $\psi(\alpha_2^{p,q})$ are linearly
independent---that $\psi(\beta_{p,q})$ and $\psi(\gamma_{p,q})$ are
nonzero can be proved using a similar argument. 

Suppose $q$ is obtained from $p$ by adding three distinct parts of
$p$. For $i\in\{1,2\}$, let $P_i=(X\arrow Y_i\arrow Z)$ be the paths
used to define $\psi(\alpha_i^{p,q})=\iNorm{P_i}$ above. If
$\lambda_1\psi(\alpha_{1}^{p,q})=\lambda_2\psi(\alpha_{2}^{p,q})$,
then $\lambda_1\Norm(P_1)-\lambda_2\Norm(P_2)$ is an element of
$\ker(\varphi)$. Thus, $Z(\lambda_1\Norm(P_1)-\lambda_2\Norm(P_2))X\in
Z(\ker\varphi)X$. By Theorem \ref{t:WEquivariantPhi},
$Z(\ker\varphi)X$ is spanned by $P_1+P_2+P_3$, where $P_3$ is the
third path from $X$ to $Z$. Hence, either $P_3$ is in the orbit of
$P_1$ or $P_2$, or $\Norm(P_3)=0$. The latter happens if and only if
$|B|=|C|$, contradicting that $|B|=p_2\neq p_3=|C|$. The
former happens if and only if $|A|=|B|$ or $|A|=|C|$. This is
again a contradiction. So
$\lambda_1\psi(\alpha_1^{p,q})\neq\lambda_2\psi(\alpha_2^{p,q})$.

Therefore, if $a\in\ker(\psi)$, then $a$ is a linear combination of
paths of length at least two. So $\ker(\psi)$ is an admissible
ideal of $k\InvariantQuiver$.
\end{proof}

\subsection{Descending Loewy Series of $\InvariantSubalgebraB$}

Here we prove the following result on the square of the radical of
$\InvariantSubalgebraB$ that was used in the proof of Theorem
\ref{t:QuiverOfDesAlgB}. The proof of this result can be adapted to
prove the corresponding result in type $A$. See
\S\ref{ss:LoewySeriesTypeA} for more details.

\begin{Theorem}
\label{t:TypeBLoewySeries}
Let $k\Faces$ be the face semigroup algebra of the reflection
arrangement of $\Bn$. Then for all $m\in\mb N$,
\begin{align*}
\rad^m\left(\InvariantSubalgebraB\right) =
\rad^{2m}\left(k\Faces\right) \cap \InvariantSubalgebraB.
\end{align*}
\end{Theorem}

During the course of the proof we will need bases of the subspaces
$X\in\IntersectionLattice$. We use a basis described by the set
partition $\pi(X) = \SetPartitionB$: for $i\in[r]$, let 
\begin{align*}
\bs\beta_{B_i} = \sum_{j\in B_i} \vec e_j,
\end{align*}
where $\vec e_1,\ldots,\vec e_n$ is the standard basis of $\mb R^n$
and $\vec e_{-j} = -\vec e_j$ for $j\in[n]$. We call 
$\bs\beta_{B_1}$, $\bs\beta_{B_2}$, \ldots, $\bs\beta_{B_r}$ 
the \defn{standard basis} of $X$.

We begin with the following lemma, which we will use several times in
the proof.

\begin{Lemma}
\label{l:ExpansionLemmaTypeB}
If $\Partition1$ is obtained from $\Partition0$ by merging two
nonzero blocks $A$ and $B$ and if $\Partition2$ is obtained from
$\Partition1$ by merging $A\cup B$ with a nonzero block $C$, then
\begin{align*}
\lambda \Norm\Path Xm & = 
\Norm(X_2\arrows X_m)\Norm(X_0\arrow X_1\arrow X_2) \\
&\qquad- \smashoperator{\sum_{t\in \Bn\atop {t(X_2)=X_2\atop t(A\cup
B\cup C) \neq \pm(A\cup B\cup C)}}}
\sigma_{X_2}(t)\sigma_{X_m}(t)\Norm(X_0\arrow X_1 \arrow t(X_2) \arrows t(X_m)),
\end{align*}
where $\lambda$ is the cardinality of
$\{t\in \Bn: t(X_2)=X_2 \text{ and } t(A\cup B\cup C) = \pm(A\cup
B\cup C)\}$.
\end{Lemma}

\begin{proof}
Note that
\begin{align*}
\Norm(X_2\arrows X_m)&\Norm(X_0\arrow X_1\arrow X_2) 
=\Norm\Big(\Norm(X_2\arrows X_m)(X_0\arrow X_1\arrow X_2)\Big) \\
&= \sum_{t\in \Bn\atop t(X_2)=X_2}
\sigma_{X_2}(t)\sigma_{X_m}(t)\Norm(X_0\arrow X_1 \arrow t(X_2) \arrows t(X_m)).
\end{align*}
Therefore, we need only show that if $t(X_2)=X_2$ and $t(A \cup B\cup
C) = \pm(A \cup B\cup C)$, then the summand in the above
sum is $\Norm\Path Xm$.

Suppose $t(X_2)=X_2$ and suppose that $t(A \cup B \cup C) =
\varepsilon(A \cup B \cup C)$, where $\varepsilon=\pm1$. Let $s$ be
the signed permutation defined by
\begin{align*}
s(i) =
\begin{cases}
\varepsilon i, & 
 \text{if } i\in A \cup B \cup C \cup \overline A \cup \overline B
 \cup \overline C, \\
t(i), & \text{otherwise.}
\end{cases}
\end{align*}
Then $s(A)=\varepsilon A, s(B)=\varepsilon B$ and $s(C)=\varepsilon C$. Hence, $s(X_i)=X_i$ for
$i\in\{0,1,2\}$ and $s(X_j)=t(X_j)$ for all $j\in\{2,\dots,m\}$. 

Next we argue that
$\sigma_{X_2}(t)\sigma_{X_m}(t)=\sigma_{X_0}(s)\sigma_{X_m}(s)$. Let
$\bs\beta_A, \bs\beta_B,\bs\beta_C, \bs\beta_1,\ldots,\bs\beta_r$
denote the standard basis for $X_0$. Then the standard basis for
$X_2$ is $\bs\beta_A+\bs\beta_B+\bs\beta_C,
\bs\beta_1,\ldots,\bs\beta_r$. Since both $s$ and $t$ induce the same
permutation on the standard basis vectors of $X_2$ and $X_m$, it
follows that $\sigma_{X_2}(s)=\sigma_{X_2}(t)$ and
$\sigma_{X_m}(s)=\sigma_{X_m}(t)$. And since $s$ either fixes or
negates all three vectors $\bs\beta_A,\bs\beta_B,\bs\beta_C$ it
follows that $\sigma_{X_0}(s)=\sigma_{X_2}(s)$. Thus,
\begin{align*}
\sigma_{X_2}(t)\sigma_{X_m}(t)&\Norm(X_0\arrow X_1 \arrow t(X_2)
\arrows t(X_m))\\
&= \sigma_{X_0}(s)\sigma_{X_m}(s)\Norm(s(X_0)\arrow s(X_1) \arrow s(X_2)
\arrows s(X_m))\\
&= \Norm(s\Path Xm) = \Norm\Path Xm.\qedhere
\end{align*}
\end{proof}

\begin{proof}[Proof of Theorem \ref{t:TypeBLoewySeries}]
\def\lhs{\ensuremath{\rad^2\left(\InvariantSubalgebra[\Bn]\right)}}
\def\rhs{\ensuremath{\rad^4\left(k\Faces\right) \cap
\InvariantSubalgebra[\Bn]}}

We first argue that we need only prove the cases $m=1,2$. 

\textbf{Reduction to the cases $\bs{m=1,2}$.}
Let $\psi:k\InvariantQuiver \to \InvariantSubalgebra[\Bn]$ be the
$k$-algebra homomorphism defined in the proof of Theorem
\ref{t:QuiverOfDesAlgB}. Note that $\psi=\varphi\circ\xi$, where $\xi$
maps paths of length $l$ in $\InvariantQuiver$ to paths of length $2l$
in $\Quiver$. The case $m=2$ was what was needed to prove that $\psi$
is surjective. Hence, if $a \in
\rad^{2m}(k\Faces)\cap\InvariantSubalgebra[\Bn]$ for some $m\in\mathbb
N$, then there exists $c\in k\InvariantQuiver$ such that $\psi(c)=a$.
We will argue that $c\in\rad^m(k\InvariantQuiver)$, thus showing that
$a=\psi(c)\in\psi(\rad^m(k\InvariantQuiver))
\subseteq\rad^m\InvariantSubalgebra[\Bn]$. 
(The reverse containment is immediate.)

Since $\varphi(\rad^p(k\Quiver))=\rad^p(k\Faces)$ for all $p\in\mb N$
(this follows from the fact that $\Quiver$ is the quiver of $k\Faces$
and contains no oriented cycles \cite[Corollary II.2.11]{Assem2006}),
it follows that $\xi(c)$ is a linear combination of paths of $\Quiver$
having length at least $2m$. Hence, $c$ is a linear combination of
paths of $\InvariantQuiver$ having length at least $m$, 
so $c \in \rad^m(k\InvariantQuiver)$.

\textbf{The case $\bs{m=1}$.}
In the proof of Proposition \ref{p:MyAlgebrasAreBasic} we argued that
$\rad(\InvariantSubalgebra)=\rad(k\Faces)\cap\InvariantSubalgebra$
for any finite Coxeter group $W$. So we need only show that
$\rad(k\Faces)\cap\InvariantSubalgebra[\Bn]=\rad^2(k\Faces)\cap\InvariantSubalgebra[\Bn]$.
Let $\varphi:k\Quiver\to k\Faces$ denote the surjection of Theorem
\ref{t:WEquivariantPhi}. Then $\rad^i(k\Faces)$ is spanned by the
elements $\varphi(P)$, where $P$ is a path of length at least $i$.
Since the transformation $\vec v\mapsto -\vec v$ is an element of
$\Bn$, it follows that $\sum_{w\in \Bn} w(P)=0$ if $P$ is a path of
odd length (see Proposition \ref{p:CoverRelationsNotGivingArrows} and
its proof). So
$\rad^{2i}(k\Faces)\cap\InvariantSubalgebra[\Bn]=\rad^{2i-1}(k\Faces)\cap\InvariantSubalgebra[\Bn]$
for $i\geq1$ since both are spanned by the elements
$\varphi(\sum_{w\in \Bn}w(P))$, where $P$ is a path of length at least
$2i$.

\textbf{The case $\bs{m=2}$.}
We first argue that $\lhs \subseteq \rhs$. Let $a\in\lhs$. Then $a=bc$
for two elements $b,c\in\rad\InvariantSubalgebra[\Bn]
=\rad^2(k\Faces)\cap\InvariantSubalgebra[\Bn]$. Thus, $bc$ is an
element of $\rad^4(k\Faces)$ and $\InvariantSubalgebra[\Bn]$. 

We prove the reverse containment by contradiction. Suppose $\lhs
\subsetneq \rhs$. Since $\rhs$ is spanned by elements of the form
$\iNorm{P}$, where $P$ is a path in $\Quiver$ of length at least
$m\geq4$, it follows that there exists a path $P=\Path Xm$ such that
$m\geq 4$ and $\iNorm P\notin\lhs$. 
We first argue that we can assume that $P$ satisfies the following:
$\Partition1$ contains a nonzero block $A\cup B$, where $A\neq B$ are
blocks of $\Partition0$ (Step 1);
$\Partition2$ contains the nonzero block $A\cup B\cup C$, where
$C$ is a block of $\Partition1$ (Step 2);
$\Partition3$ contains a nonzero block $D\cup E$, where $D$ and $E$
are blocks of $\Partition2$ and $|D|=|A\cup B\cup C|$ (Step 3).
Then we derive a contradiction (Step 4).

\textbf{Step 1.} We argue that $\Partition1$ contains a nonzero block
$A\cup B$, where $A\neq B$ are blocks of $\Partition0$.

If not, then the zero block $Z_1$ of $\Partition1$ is $Z_1 = B \cup
Z_0 \cup \overline B$ for some nonzero block $B\in\Partition0$, where
$Z_0$ is the zero block of $\Partition0$. Let $t$ be the signed
permutation that negates the elements of $B$ and fixes the other
elements. Then $\Norm(P) = \Norm(t(P)) = -\Norm(P)$. Hence,
$\Norm(P)=0$, contradicting that $\iNorm P\notin\lhs$.

\textbf{Step 2.}
We argue that $\Partition2$ contains the nonzero block $A\cup B\cup
C$, where $C$ is a block of $\Partition1$.

First we show that $A \cup B$ is not a block of $\Partition m$. If $A
\cup B$ is a block of $\Partition m$, then it is a block of
$\Partition j$ for all $j\in\{1,\ldots,m\}$. Let $t$ be the signed
permutation that negates the elements of $A\cup B$. Then $t(P)=-P$, so
$\Norm(P)=0$, contradicting that $\iNorm P\notin\lhs$. 

This implies that there exists $j\in[m]$ such that $\Partition j$ is
obtained from $\Partition{j-1}$ by merging $A \cup B$. We argue that
we can assume $j=2$. From Theorem \ref{t:WEquivariantPhi}, it follows
that \begin{align*} P \quad+ \sum_{Y\neq X_{j-1} \atop X_{j-2}\arrow Y
\arrow X_j} \big(X_0 \arrows X_{j-2} \arrow Y \arrow X_j \arrows
X_m\big) \in \ker(\varphi). \end{align*} Since $\iNorm P\notin\lhs$,
there must exist at least one $Y\neq X_{j-1}$ such that $\iNorm{X_0
\arrows X_{j-2} \arrow Y \arrow X_j \arrows X_m}\notin\lhs$. By
examining the intervals of length two in $\PartitionLatticeB$ (see
Figure \ref{f:PartitionLatticeTypeBTwoIntervals}), we note that
$\pi(Y)$ is obtained from $\Partition{j-2}$ by merging $A \cup B$ with
some other block of $\Partition{j-2}$. By replacing $P$ with this
path, noticing that this new path still begins with $X_0\arrow X_1$,
and repeating this argument until $j=2$, we have that $\Partition2$ is
obtained from $\Partition1$ by merging $A\cup B$.

If the zero block of $\Partition2$ is $Z_2 = A\cup B \cup Z_0 \cup
\overline A \cup \overline B$, where $Z_0$ is the zero block of
$\Partition0$, then any $t\in\Bn$ that fixes $X_2$ also fixes $Z_2$.
This implies, by appealing to the argument in the proof of Lemma
\ref{l:ExpansionLemmaTypeB}, that $\Norm(P)$ is a scalar multiple of
$\Norm(X_2\arrows X_m)\Norm(X_0\arrow X_1\arrow X_2)$, contradicting
that $\iNorm P\notin\lhs$. 

\textbf{Step 3.} We argue that we can suppose that $\Partition3$
contains a nonzero block $D\cup E$, where $D$ and $E$ are blocks of
$\Partition2$ and $|D|=|A\cup B\cup C|$.

We first argue that not all of the blocks of $\Partition2$ of size
$\lambda = |A \cup B \cup C|$ are also blocks of $\Partition m$.
We do this by
showing that we can factor $\Norm(P)$. By Lemma
\ref{l:ExpansionLemmaTypeB}, we need only show that if $t(X_2)=X_2$,
then 
\begin{gather*}
\sigma_{X_2}(t)\sigma_{X_m}(t)\Norm\Big(X_0\arrow X_1\arrow t(X_2) \arrows t(X_m)\Big) = \Norm(P).
\end{gather*}

Write $\Partition2 = \{B_1,\ldots,B_k,C_1,\ldots,C_l; Z_0\},$ where
$|C_i|=\lambda$, $|B_j|\neq\lambda$, and write $\Partition m =
\{D_1,\ldots,D_h,C_1,\ldots,C_l; Z_m\}.$ Suppose $t\in \Bn$ such that
$t(X_2)=X_2$. Then
$t$ permutes the blocks $\pm B_1,\ldots,\pm B_k$, as well as the
blocks $\pm C_1,\ldots,\pm C_l$. Define $s\in \Bn$ by
$s|_{B_i}=t,s|_{Z_0}=t$ and $s|_{C_j}=1$. Then $s(X_0)=X_0$,
$s(X_1)=X_1$ and $s(X_j)=t(X_j)$ for all $j\in\{2,\ldots,m\}$. It
remains to show that $\sigma_{X_0}(s)\sigma_{X_m}(s)=
\sigma_{X_2}(t)\sigma_{X_m}(t)$. This follows by
comparing the actions of $s$ and $t$ on the standard basis
$\bs\beta_{B_1},\ldots,\bs\beta_{B_k},\bs\beta_{C_1},\ldots,\bs\beta_{C_l}$
of $X_2$.

This implies that some nonzero block $D$ of $\Partition1$ of size
$\lambda$ is merged to get $\Partition j$ for some
$j\in\{3,\ldots,m\}$. If the zero block of $\Partition j$ is $Z_j=D\cup
Z_{j-1}\cup \overline{D}$, where $Z_{j-1}$ is the zero block of
$\Partition{j-1}$, the let $t$ be the signed permutation that negates
the elements of $D$ and fixes the other elements of $[\pm n]$. Note
that since $D$ is a block of $\Partition1$ and $|D|=|A\cup B\cup C|$,
we have either that $D$ is a block of $\Partition0$ or $D=A\cup B\cup
C$. In both cases $t$ negates an odd number of elements of the
standard basis of $X_0$ and no elements of the standard basis of
$X_m$. Thus, $\sigma_{X_0}(t)=-1$ and $\sigma_{X_m}(t)=1$. It follows
that $t(P)=-P$, and so $\iNorm P = 0$, a contradiction.

By arguing as in Step 2, using the relations in $\ker(\varphi)$, we
can assume that $j=3$. 

\textbf{Step 4.} We are now ready to conclude the proof. Let $P=\Path
Xm$ be a path of length $m\geq4$ such that $\iNorm P\notin\lhs$. From
Step 1 we have that $\Partition1$ contains a nonzero block $A \cup B$,
where $A\neq B$ are blocks of $\Partition0$. Of all such paths, pick
$P$ such that $|A|+|B|$ is maximal. That is, we suppose the following.
\begin{enumerate}
\item[($\star$)]
If $P'=(X_0\arrow Y_1 \arrow Y_2 \arrows Y_m)$ with
$\iNorm{P'}\notin\lhs$ and if $\pi(Y_1)$ contains
the nonzero block $A'\cup B'$, where $A'\neq B'$
are blocks of $\Partition0$, then $|A'\cup B'| \leq |A\cup B|$.
\end{enumerate}
In both Steps 2 and 3, we replaced $P$ with other paths $P'$ that
begin with $X_0\arrow X_1$ and such that $\iNorm{P'}\notin\lhs$.
Therefore, we can assume that $\Partition2$ contains the nonzero block
$A\cup B\cup C$, where $C$ is a block of $\Partition1$ and that
$\Partition3$ contains a nonzero block $D\cup E$, where $D$ and $E$
are blocks of $\Partition2$ and $D$ has cardinality $\lambda=|A\cup
B\cup C|$. Therefore, there are three cases to consider.

\def\CaseOneFig{%
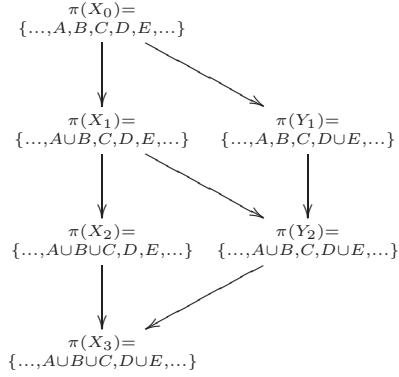
\begin{figure}[!htb]
\ifhrule\hrule\fi
{\footnotesize 
\begin{gather*}
\xymatrix@C=0em{
{\Partition0\mathord=\atop\{\cldots,A,B,C,D,E,\cldots\} }
 \ar[d] \ar[dr] \\
{\Partition1\mathord=\atop\{\cldots,A \cup B, C, D, E, \cldots\} }
 \ar[d] \ar[dr] 
&{\pi(Y_1)\mathord=\atop\{\cldots,A,B, C, D \cup E,\cldots\}}
 \ar[d] \\
{\Partition2\mathord=\atop\{\cldots,A \cup B \cup C, D, E,\cldots\} }
 \ar[d] 
&{\pi(Y_2)\mathord=\atop\{\cldots,A \cup B, C, D \cup E,\cldots\} }
 \ar[dl] \\
{\Partition3\mathord=\atop\{\cldots,A \cup B \cup C, D \cup E,\cldots\}}
}
\end{gather*}}
\caption{In Case 1, $D$ and $E$ are not $\pm(A\cup B \cup C)$.
Note that $|D \cup E| > |A\cup B|$.}
\label{f:case1}
\ifhrule\hrule\fi
\end{figure}}
\def\CaseTwoFig{%
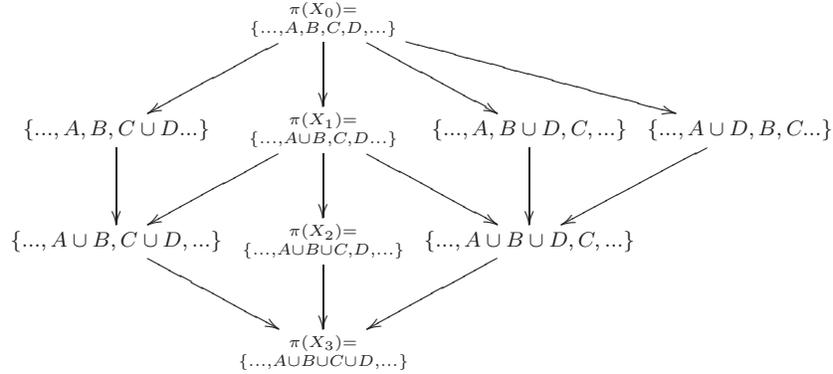
\begin{figure}[!htb]
\ifhrule\hrule\fi
{\footnotesize 
\begin{gather*}
\xymatrix@C=0em{
&& {\Partition0=\atop\{\cldots,A,B,C,D,\cldots\}}
 \ar[d] \ar[dl] \ar[dr] \ar[drr] \\
& \{\cldots,A, B, C \cup D \cldots\} 
 \ar[d]
& {\Partition1=\atop\{\cldots,A \cup B, C, D \cldots\}}
 \ar[d] \ar[dr] \ar[dl] 
& \{\cldots,A, B \cup D, C, \cldots\} 
 \ar[d]
& \{\cldots,A \cup D, B, C \cldots\} 
 \ar[dl]
\\
& \{\cldots,A \cup B, C \cup D, \cldots\} 
 \ar[dr]
& {\Partition2=\atop\{\cldots,A \cup B \cup C, D,\cldots\}}
 \ar[d] 
& \{\cldots,A \cup B \cup D, C,\cldots\} 
 \ar[dl] \\
&& {\Partition3=\atop\{\cldots,A \cup B \cup C \cup D,\cldots\}}
}
\end{gather*}}
\caption{In Case 2, $|D|=\lambda$, thus
$|C\cup D|, |B\cup D|, |A \cup D| > |A \cup B|$.}
\label{f:case2}
\ifhrule\hrule\fi
\end{figure}}
\def\CaseThreeFig{%
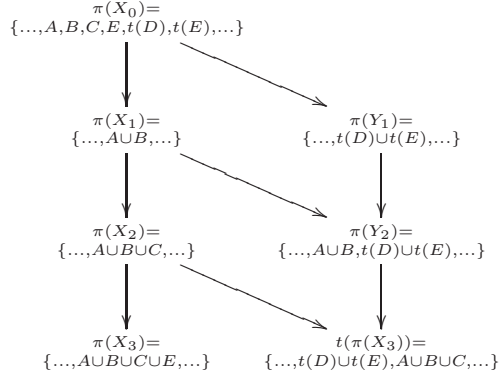
\begin{figure}[!htb]
\ifhrule\hrule\fi
{\footnotesize 
\begin{gather*}
\xymatrix@C=0em{
{\Partition0\mathord=\atop\{\cldots,A,B,C,E,t(D),t(E), \cldots\} }
 \ar[d] \ar[dr] \\
{\Partition1\mathord=\atop\{\cldots,A \cup B, \cldots\} }
 \ar[d] \ar[dr] 
&{\pi(Y_1)\mathord=\atop\{\cldots,t(D) \cup t(E), \cldots\}}
 \ar[d] \\
{\Partition2\mathord=\atop\{\cldots,A\cup B\cup C,\cldots\}}
 \ar[dr] 
&{\pi(Y_2)\mathord=\atop\{\cldots, A\cup B, t(D)\cup t(E), \cldots\} }
\ar[d] \\
{\Partition3\mathord=\atop\{\cldots,A\cup B\cup C\cup E,\cldots\}}\ar@{<-}[u]
& {t(\Partition3)\mathord=\atop\{\cldots, t(D)\cup t(E), A\cup B \cup C, \cldots\}}
}
\end{gather*}}
\caption{In Case 3, $|E|\neq\lambda$. 
Note that $|t(D)\cup t(E)|>|A \cup B|$.}
\label{f:case3}
\ifhrule\hrule\fi
\end{figure}}

\emph{Case 1.} Suppose $D,E\neq \pm(A\cup B\cup C)$. This case is
illustrated in Figure \ref{f:case1}. \CaseOneFig The open
interval $(X_3,X_1)=\{Y\in\IntersectionLattice:X_3<Y<X_1\}$ contains
exactly two elements: $X_2$ and $Y_2$, where $\pi(Y_2)$ is obtained
from $\Partition1$ by merging $D$ with $E$. The open interval
$(Y_2,X_0)$ also contains exactly two elements: $X_1$ and $Y_1$, where
$\pi(Y_1)$ is obtained from $\Partition0$ by merging $D$ and $E$. If
$P'=(X_0\arrow Y_1 \arrow Y_2 \arrow X_3 \arrows X_m)$, then
$\iNorm{P'}\notin\lhs$ since $P'-P \in\ker\varphi$ (Lemma
\ref{l:KernelOfPhi}), but $|D\cup
E|>|A\cup B|$, contradicting ($\star$).

\emph{Case 2.} Suppose $D\neq E=A\cup B\cup C$. This situation is
illustrated in Figure \ref{f:case2}. \CaseTwoFig 
Let $(X_0\arrow Y_i\arrow Z_i\arrow X_3)$ for $i\in\{1,2,3\}$ be the
three paths in Figure \ref{f:case2} from $X_0$ to $X_3$ such that
$Y_i\neq X_1$ and $Z_i\neq X_2$, and let $P_i = (X_0\arrow Y_i\arrow
Z_i\arrow X_3\arrows X_m)$ for $i\in\{1,2,3\}$. Since $|D|>|A\cup B|$,
the assumption ($\star$) implies that $\iNorm{P_i}\in\lhs$ for
$i\in\{1,2,3\}$. But $P-(P_1+P_2+P_3)\in\ker(\varphi)$ (Lemma
\ref{l:KernelOfPhi}), so this
implies that $\Norm(P)\in\lhs$, a contradiction.

\emph{Case 3.} Suppose $D=A\cup B\cup C\neq E$. If $|E|=\lambda$,
then we can swap the roles of $D$ and $E$ and apply the argument
from Case 2. So suppose that $|E|\neq\lambda$. By Lemma
\ref{l:ExpansionLemmaTypeB}, there exists $t\in \Bn$ such that
$t(X_2)=X_2$, $t(D)\neq\pm D$ and 
\begin{gather}
\label{e:Case3AssumptionTypeB}
\iNorm{X_0\arrow X_1
\arrow X_2 \arrow t(X_3) \arrows t(X_m)}\notin\lhs.
\end{gather}

We argue that we are in the situation illustrated in Figure
\ref{f:case3}.
\CaseThreeFig
We first establish that $t(D), t(E) \in \Partition0$.
Since $t(X_2) = X_2$, it follows that $t$ permutes the blocks of
$\Partition2$. Since $t(D)\neq\pm D$, it follows that $t(D)$ is a
block of $\Partition2$ different than $\pm D = \pm(A\cup B \cup C)$.
And because all other blocks of $\Partition2$ are blocks of
$\Partition0$, we have that $t(D)$ is a block of $\Partition0$.
Considering that $|E|\neq\lambda = |D|$, we have $t(E)\neq\pm D$, so
the same reasoning implies that $t(E)$ is also a block of
$\Partition0$.

Since $t(X_2) = X_2$ and $\Partition3\lessdot\Partition2$, it follows
that $t(\Partition3)\lessdot\Partition2$. Therefore, $t(\Partition3)$
is obtained from $\Partition2$ by merging $t(D)$ and $t(E)$ since
$t(D),t(E)\in\Partition2$ and $t(D)\cup t(E) = t(A\cup B \cup C \cup
E) \in t(\Partition3)$. There is exactly one other partition
$\pi(Y_2)$ such that
$t(\Partition3)\lessdot\pi(Y_2)\lessdot\Partition1$, the partition
obtained from $\Partition1$ by merging $t(D)$ with $t(E)$. There is
exactly one other partition $\pi(Y_1)$ such that
$\pi(Y_2)\lessdot\pi(Y_1)\lessdot\Partition0$, the partition obtained
from $\Partition0$ by merging $t(D)$ with $t(E)$. So we are in the
situation illustrated in the figure.

By Lemma \ref{l:KernelOfPhi}, the following element is in
$\ker(\varphi)$:
\begin{align*}
(X_0\arrow X_1 \arrow X_2 \arrow t(X_3) \arrows t(X_m))
-(X_0\arrow Y_1 \arrow Y_2 \arrow t(X_3) \arrows t(X_m)).
\end{align*}
Together with \eqref{e:Case3AssumptionTypeB}, this implies that 
\begin{gather*}
\iNorm{X_0\arrow Y_1 \arrow Y_2 \arrow t(X_3) \arrows t(X_m)}\notin\lhs.
\end{gather*}
This contradicts our assumption ($\star$) because $\pi(Y_1)$ is
obtained from $\Partition0$ by merging the blocks $t(D)$ and $t(E)$,
and $|t(D)\cup t(E)| = |t(A\cup B \cup C \cup E)| > |A \cup B|$.
\end{proof}

\section{Future Directions}
\label{s:FutureDirections}

This article is part of an ongoing project to determine the quiver
with relations of the descent algebras. There is still much to do.

The quivers of all the descent algebras have not yet been determined.
The main outstanding case is the quiver of the descent algebra of type
$D$ as the exceptional types can be dealt with using computer algebra
software \cite{Pfeiffer2007}. It should be possible to adapt the
proofs of Theorems \ref{t:QuiverOfDesAlgA} and \ref{t:QuiverOfDesAlgB}
to this case as well, but given the similarity between these
arguments, a general argument is more desirable. The main obstacle 
is to understand the relationship between
$\rad^2(\InvariantSubalgebra)$ and
$\InvariantSubalgebra\cap\rad^p(k\Faces)$. Indeed, if these two spaces
are equal for some $p$, then there is no arrow from $\Orbit'$ to
$\Orbit$, where $\Orbit', \Orbit\in\IntersectionLattice/W$, if
$\Orbit<\Orbit'$ and $\rank(\Orbit')-\rank(\Orbit)\geq p$; so only the
intervals in $\IntersectionLattice/W$ of length $p-1$ need to be
studied. This is precisely what we did for types $A$ and $B$, where
$p$ was $2$ and $4$, respectively.

Another task is to determine relations for some quiver presentation of
the descent algebras. Very little is known here, even for type $A$.

Other representation theoretic questions also arise. As mentioned
following Corollary \ref{c:DescentAlgebraIsQuasiHereditary}, it would
be interesting to determine the characteristic tilting module of each
descent algebra as well as its Ringel dual. Also, the Cartan
invariants of the descent algebras are not known in general. Formulas
exist for type $A$ (see \cite{GarsiaReutenauer1989},
\cite[Corollary~2.1]{BlessenohlLaue1996},
\cite[Section~3.6]{KrobLeclercThibon1997} and
\cite[Section~9.4]{Schocker2006}) and a combinatorial interpretation
for type $B$ was given by Nantel Bergeron
\cite[Theorem~3.3]{NBergeron1992}.

%\bibliographystyle{halpha}
%\bibliography{references} 

\end{document}